\documentclass[preprint,10pt,final]{elsarticle}
\usepackage[letterpaper]{geometry}
\usepackage[usenames,dvipsnames]{color}
\usepackage{bm}
\usepackage{soul}
\usepackage[ampersand]{easylist}
\usepackage{amsmath}
\usepackage{amssymb}
\usepackage{amsfonts}
\usepackage{graphicx}
\usepackage{parskip}
\usepackage{longtable}
\usepackage{booktabs}
\usepackage{subcaption}
\usepackage[usenames,dvipsnames]{color}
\usepackage{multirow}
\usepackage[mathscr]{euscript}
\usepackage{algorithmicx}
    \usepackage[ruled]{algorithm}
    \usepackage{algpseudocode}
\biboptions{sort&compress}

\newcommand{\mcal}[1]{\mathcal{#1}}

\newcommand{\sbf}[1]{\boldsymbol{#1}}
\newcommand{\abs}[1]{\bigl\lvert{#1}\bigr\rvert}

\newcommand{\lp}{\left(}
\newcommand{\rp}{\right)}
\newcommand{\blp}{\bigl(}
\newcommand{\brp}{\bigr)}
\newcommand{\lb}{\left[}
\newcommand{\rb}{\right]}

\newcommand{\bx}{\boldsymbol{x}}

\newcommand{\bn}{\boldsymbol{n}}

\newcommand{\beql}{\begin{equation}}
\newcommand{\eeql}{\end{equation}}
\newcommand{\beq}{\begin{equation*}}
\newcommand{\eeq}{\end{equation*}}
\newcommand{\bali}{\begin{align}}
\newcommand{\eali}{\end{align}}
\newcommand{\balit}{\begin{alignat}}
\newcommand{\ealit}{\end{alignat}}
\newcommand{\bse}{\begin{subequations}}
\newcommand{\ese}{\end{subequations}}
\newcommand{\bfig}{\begin{figure}[!ht]}
\newcommand{\efig}{\end{figure}}

\newcommand{\eMat}{\end{matrix} \right]}

\providecommand{\apo}{\emph{a posteriori }}
\providecommand{\apr}{\emph{a priori }}

\newcommand{\ggrad}{\boldsymbol{\nabla}}

\def\etal{{\it et al. }}

\journal{Journal of Computational Physics}

\begin{document}

\begin{frontmatter}

\title{A model and variance reduction method for computing statistical outputs of stochastic elliptic partial differential equations}

\author[mit]{F.~Vidal-Codina\corref{cor1}\fnref{fn1,fn2}}
\ead{fvidal@mit.edu}
\author[mit]{N.~C.~Nguyen\fnref{fn1}}
\ead{cuongng@mit.edu}
\author[ox]{M.~B.~Giles}
\ead{mike.giles@maths.ox.ac.uk}
\author[mit]{J.~Peraire\fnref{fn1}}
\ead{peraire@mit.edu}

\cortext[cor1]{Corresponding author}
 \address[mit]{Department of Aeronautics and Astronautics, Massachusetts Institute of Technology, Cambridge, MA 02139, USA}
\address[ox]{Mathematical Institute, University of Oxford, Oxford, UK}
\fntext[fn1]{This work was supported by AFOSR Grant No. FA9550-11-1-0141, AFOSR Grant No. FA9550-12-0357, and the Singapore-MIT Alliance.}
\fntext[fn2]{This work was supported by Obra Social la Caixa.}

 \begin{abstract}
We present a model and variance reduction method for the fast and reliable computation of statistical outputs of stochastic elliptic partial differential equations. Our method consists of three main ingredients: (1) the hybridizable discontinuous Galerkin (HDG) discretization of elliptic partial differential equations (PDEs), which allows us to obtain high-order accurate solutions of the governing PDE; (2) the reduced basis method for a new HDG discretization of the underlying PDE to enable real-time solution of the parameterized PDE in the presence of stochastic parameters; and (3) a multilevel variance reduction method that exploits the statistical correlation among the different reduced basis approximations and the high-fidelity HDG discretization to accelerate the convergence of the Monte Carlo simulations. The multilevel variance reduction method provides  efficient computation of the statistical outputs by shifting most of the computational burden from the high-fidelity HDG approximation to the reduced basis approximations.  Furthermore, we develop \textit{a posteriori} error estimates for our approximations of the statistical outputs. Based on these error estimates, we propose an algorithm for optimally choosing both the dimensions of the reduced basis approximations and the sizes of Monte Carlo samples to achieve a given error tolerance.  We provide numerical examples to demonstrate the performance of the proposed method.
 \end{abstract}
  
\end{frontmatter}

\section{Introduction}

The analysis of physical systems is often carried out by mathematical modeling and numerical simulation. For a given system, the corresponding mathematical model requires certain input data, such as material properties, forcing terms, boundary conditions and geometry information. For many problems of interest, input data are not known precisely. In such cases, one may need to consider input data as random variables and represent them in probabilistic terms. Mathematical models represented by partial differential equations with random input data are known as stochastic partial differential equations (SPDEs). Uncertainty in the input data may come from different sources. It can be that the physical system has some intrinsic variability, for example, uncertainty in the gust loads on an aircraft, or wind and seismic loading on civil structures. It is also possible that we are unable to effectively characterize the physical system with a mathematical model because, for instance, we 
may have errors in geometry, roughness of surfaces, or multiscale behavior that we are unable to capture.  Therefore, there is a growing need to represent the uncertainty in the data and effectively propagate it through the mathematical model. The goal of this probabilistic approach resides in computing statistics of some observable outputs (quantities of interest), which are usually defined as functionals of the solution of the underlying SPDE.

There exist a number of different approaches to solve SPDEs and retrieve the statistics of the output. The most common approach is to use Monte Carlo (MC) methods \cite{fishman1996,hammersley1965monte}. Monte Carlo methods only need repeated evaluations of the output functional of the solution of the SPDEs for different instantiations of the random input. The main advantage of Monte Carlo methods is that their convergence rate is independent of the dimension of the stochastic space, namely, the number of random variables. The main caveat of these methods is their slow convergence rate, which demands a large amount of realizations to achieve accurate results. As a result, a number of techniques such as quasi Monte Carlo methods \cite{nieder1992,caflisch1998monte}, Latin Hypercube Sampling \cite{loh1996LHS,stein1987LHS}, variance reduction methods  \cite{caflisch1998monte} and multilevel Monte Carlo \cite{giles2008multilevel} have been proposed to alleviate the slow convergence rate of the standard Monte Carlo 
methods. 

Another approach is stochastic Galerkin methods, first introduced by Ghanem \etal in \cite{ghanem1991}, that generalize the theory of Wiener-Hermite polynomial chaos expansion \cite{wiener} and combine it with a finite element method to model uncertainty in a SPDE. In this approach, the random variables are treated as additional dimensions of the problem and projected onto a stochastic space spanned by a set of orthogonal polynomials. The problem is then reduced to a system of deterministic equations, which couple the physical and stochastic dimensions. This methodology has proven to be very effective when solving partial differential equations (PDEs) in a broad range of applications, such as diffusion problems and heat conduction \cite{ghanem1999,xiu2002modeling,xiu2003new}, structural dynamics \cite{sarkar2002}, transport in random media \cite{ghanemheterogeneous} and fluid dynamics \cite{xiu2003modeling,chen2005steady}. The advantage of these methods is that they converge exponentially fast for a sufficiently regular solution 
field \cite{babuska2004galerkin,babuvska2005solving,deb2001solution}. However, their main drawback is that their computational complexity  grows combinatorially with the number of random variables and the number of expansion terms. As a consequence, they are not effective  for solving problems with a large number of random variables.

A more recent approach is stochastic collocation methods (SCM), first introduced in \cite{mathelin2003} and further developed in \cite{xiu2005hocollocation}. The main idea is to compute deterministic solutions of the SPDE for certain instantiations of the random variables and then construct an interpolation function to approximate the response over the stochastic space. When the interpolation procedure is performed on tensor grids, these methods suffer from the exponential growth with the dimensionality of the stochastic space. To economize the interpolation process in large dimensions, sparse grids (Smolyak \cite{smolyak}) were introduced for elliptic problems \cite{nobile2008SG,xiu2005hocollocation}, parabolic problems \cite{nobile2009analysis} and natural convection problems \cite{ganapathy2007}. In \cite{babuvska2007stochastic}, sparse grids were shown to achieve exponential convergence for problems with smooth solutions. However, like polynomial chaos expansions, sparse grids still suffer 
from the curse of dimensionality in the sense that the number of grid points grows rapidly with the dimension of the stochastic space. Recently, anisotropy and adaptivity on sparse grids \cite{grestner2003dimadaptive,klimke2006uncertainty} have been used in SCM \cite{ganapathy2007,nobile2008anisotropic} to mitigate the elevated cost in high dimensions.

In this paper, we present a model and variance reduction method for the fast and reliable computation of statistical outputs of stochastic elliptic partial differential equations. Our method consists of three main ingredients: (1) the hybridizable discontinuous Galerkin (HDG) discretization of elliptic partial differential equations (PDEs), which allows us to obtain high-order accurate solutions of the governing PDE; (2) a reduced basis method for the HDG discretization of the underlying PDE to enable real-time solution of the parameterized PDE in the presence of stochastic parameters; and (3) a multilevel variance reduction method that exploits the statistical correlation among the different reduced basis approximations and the high-fidelity HDG discretization to accelerate the convergence rate of the Monte Carlo simulations. The multilevel variance reduction method provides  efficient computation of the statistical outputs by shifting most of the computational 
burden from the high-fidelity HDG approximation to the reduced basis approximations. Although the three ingredients of our approach exist in the literature, the main contribution of this paper is to put these methodologies into a unified framework that combines all of their strengths to tackle stochastic elliptic PDEs. Another important contribution of the paper is to develop {\em a posteriori} error bounds for the estimates of the statistical outputs and to introduce an algorithm for optimally choosing the dimensions of the reduced basis approximations and the sizes of MC samples to achieve a given error tolerance.  Last but not least, we present a new HDG formulation that enables the efficient construction of reduced basis approximations for the HDG discretization of parameterized PDEs.

The HDG method was first introduced in \cite{cockburn2009unified} for elliptic problems, subsequently analyzed in \cite{cockburn2008superconvergent,cockburn2009superconvergent,cockburn2010projection}, and later extended to a wide variety of PDEs \cite{nguyen2009linearCD,nguyen2009nonlinearCD,cockburn2011analysis,nguyen2010incns,nguyen2010stokes,nguyen2010hybridizable,nguyen2011ns,nguyen2011maxwell,nguyen2011acoustic,Nguyen2012CM,ueckermann2010high}. The HDG method is particularly effective for solving elliptic PDEs because it possesses several unique features that distinguish it from other DG methods. First, it reduces the number of globally coupled unknowns to those required to represent the trace of the approximate solution on the element boundaries, thereby resulting in a smaller global systems than other DG methods.  Second, the method provides optimal convergence rates for both the solution and the flux. And third, its flux superconvergence properties can be exploited to devise a local 
a postprocess that increases the convergence rate of the approximate solution by one order. These advantages are the main driver for the development of the Reduced Basis (RB) method for the HDG discretization of parameterized PDEs.  While the RB method is well developed for the standard finite element discretization of parameterized PDEs \cite{noor1980reduced,machiels2000output,prud2002reliable,grepl2005posteriori,grepl2005reduced,veroy2003posteriori,veroy2005certified}, the RB method for the HDG approximation of parameterized PDEs has not been considered before. The HDG discretization has multiple field variables and various equivalent weak formulations, which make the application of the RB method non straightforward.

Recently, the RB method has been applied to standard continuous Galerkin finite element solutions of stochastic elliptic PDEs \cite{boyaval2010reduced,haasdonk2013reduced,chen2014comparison}. In this approach, the stochastic PDE is first reformulated as a parametrized PDE over the coefficients of the Karhunen-Lo\`{e}ve expansion of the random fields. The reduced basis approximation and associated {\em a posteriori} error estimation are then developed for the resulting parametrized PDE. Finally, the output statistics and their error estimates are computed with a MC simulation \cite{boyaval2010reduced,haasdonk2013reduced} or a stochastic collocation approach \cite{chen2014comparison}. These approaches, which involve the RB method and its \apo error bounds to evaluate the output instead of the original finite element discretization, have been shown to outperform both standard MC and stochastic collocation. In this paper, we extend the previous work \cite{boyaval2010reduced,haasdonk2013reduced} in several important ways. We will use the HDG method to construct the RB approximation. We will adopt the multilevel Monte Carlo strategy \cite{giles2008multilevel,
barth2011multi,cliffe2011multilevel,teckentrup2013further} and demonstrate a significant computational gain relative to the standard MC approach. Moreover, we will provide {\em a posteriori} error estimates for our prediction of the statistical outputs without involving  {\em a posteriori} error bounds for the RB approximation. This feature will broaden the applicability of our approach to a wide variety of stochastic PDEs for which {\em a posteriori} error bounds for the RB approximation are either not available or too expensive to compute.

According to the central limit theorem \cite{feller1968introduction}, the error in a Monte Carlo estimation of the expectation of an output is proportional to the square root of the ratio of the variance of the output and the number of samples. Therefore, in order to reduce the error one can increase the number of samples and/or decrease the variance of the output. Because increasing the number of samples leads to higher computational cost, various techniques such as the control variates method \cite{boyaval2012fast,caflisch1998monte,hammersley1965monte}, the multilevel Monte Carlo method \cite{giles2008multilevel,barth2011multi,cliffe2011multilevel,teckentrup2013further}, and the multi-fidelity Monte Carlo method \cite{ng2012multifidelity} have been proposed to reduce the variance of the output. The control variates method reduces the variance of the output by making use of the correlation between the output and a surrogate. The multi-fidelity Monte Carlo method makes use of the statistical correlation 
between the low-fidelity (surrogate) and high-fidelity outputs to reduce the number of high-fidelity evaluations needed to achieve a given accuracy of interest. The multilevel Monte Carlo method applies the principle of control variates to a sequence of lower fidelity outputs (multigrid approximations) to estimate the statistics of the high-fidelity output. Likewise, our method applies the principle of control variates to the HDG approximation and a sequence of reduced basis approximations, thereby shifting the computational burden from the high-fidelity HDG discretization to the lower fidelity RB approximations. 

This article is organized as follows. In Section \ref{sec:RBHDG}, we introduce a stochastic elliptic boundary value problem and describe a new weak HDG formulation particularly suited for the reduced basis method. In Section 3, we describe a reduced basis method for the HDG approximation of the stochastic elliptic boundary value problem. In Section \ref{sec:RBCV}, we develop a multilevel Monte Carlo method that incorporates the HDG approximation and its reduced basis models into a unified framework to provide rapid reliable computation of the statistical outputs. In Section \ref{sec:results}, we present numerical results to demonstrate the performance of the proposed method. Finally, in Section \ref{sec:conclusions}, we discuss some directions for future research.

\section{The Hybridizable Discontinuous Galerkin Method}\label{sec:RBHDG}

\subsection{Problem statement}

Let $\mathcal{D} \in \mathbb{R}^d$ be a regular domain with Lipschitz boundary $\partial \mathcal{D}$. We consider the following stochastic boundary value problem: find a function $u$ such that,
\bse\label{eq:helm}
\begin{alignat}{2}
-\nabla\cdot (\kappa\ggrad u)  +\varrho u&= f,\qquad& \forall \bx& \in  \mathcal{D}\;, \label{eq:helmeq}\\
 \kappa \ggrad u \cdot \bn + \nu u\ &= g,\qquad& \forall \bx& \in \partial \mathcal{D} \;, \label{eq:helmR}
\end{alignat}
\ese
where $f$ is the source term, $\kappa$ is the diffusion coefficient, $\varrho$ is the Helmholtz parameter, $\nu$ is the Robin coefficient, and $g$ is the boundary data. In this problem, one or more than one of the quantities $f,\kappa,\varrho,\nu,g$ are stochastic functions. For simplicity of exposition we shall assume that $\kappa$ is a real stochastic function and that $f,\varrho,\nu,g$ are deterministic. The generalization to the case where one or more of $f,\varrho,\nu,g$ are stochastic is straightforward.   Note that since we allow $f,\varrho,\nu,g$ to be complex-valued functions, the solution $u$ is in general a complex stochastic function.

We next introduce a probability space $(\Omega,\mathcal{F},P)$, where $\Omega$ is the set of outcomes, $\mathcal{F}$ is the $\sigma$-algebra of the subsets of $\Omega$, and $P$ is the probability measure. If $Z$ is a real random variable in $(\Omega,\mathcal{F},P)$ and $\omega$ a probability event, we denote its expectation by $E[Z] = \int_{\Omega} Z(\omega) d P (\omega)$. We will consider random functions $v$ in $L^2(\mathcal{D} \times \Omega)$ equipped with the following norm
\begin{equation*}
\|v\|^2 = E\lb \int_{\mathcal{D}} |v( \bm x,\cdot)|^2 d \bm x\rb = \int_{\Omega} \int_{\mathcal{D}} |v( \bm x,\omega)|^2 d \bm x \,d P (\omega) .
\end{equation*}
We will assume that $\kappa \in L^2(\mathcal{D} \times \Omega)$ and that $\kappa( \bm x,\omega)$ is bounded and strictly positive, i.e., there exist constants $\alpha_1$ and $\alpha_2$ such that
\begin{equation*}
0 < \alpha_1 \le \kappa( \bm x,\omega) \le \alpha_2 < +\infty, \qquad \mbox{a.s. in } \mathcal{D} \times \Omega  .
\end{equation*}
We next assume that the random function $\kappa( \bm x,\omega)$ can be written in the following form
\begin{equation*}
\kappa( \bm x,\omega) = \overline{\kappa}( \bm x) + \sum_{q=1}^Q \psi_q( \bm x) y_q(\omega), \qquad ( \bm x,\omega) \in \mathcal{D} \times \Omega ,
\end{equation*}
where $\overline{\kappa}( \bm x)$ is the expectation of $\kappa$, $\psi_q( \bm x), q = 1,\ldots,Q$ are uniformly bounded real functions, and $y_q(\omega)$ for $q= 1 \ldots, Q$ are mutually independent random variables with zero mean. In addition, we assume that each of the $y_q(\omega)$ is bounded in the interval $\Lambda_q = [-\gamma_q, \gamma_q]$ with a uniformly bounded probability density function $\rho_q : \Lambda_q \to \mathbb{R}^+$. It thus follows that, with a slight overloading of notation, we can write $\kappa$ in the form 
\begin{equation*}
\label{affinekappa}
\kappa( \bm x, \bm y) = \overline{\kappa}( \bm x) + \sum_{q=1}^Q \psi_q( \bm x) y_q, \qquad ( \bm x, \bm y) \in \mathcal{D} \times \Lambda ,
\end{equation*}
where $ \bm y = (y_1,\ldots,y_Q)$ and $\Lambda = \prod_{q = 1}^Q \Lambda_q$.

Therefore, the solution $u$ of (\ref{eq:helm}) can be written as a function of $ \bm y \in \Lambda$, namely, $u( \bm x, \bm y)$. Now let $\ell$ be a bounded linear functional. We introduce a random output $s$ defined as
\begin{equation*}
s( \bm y) = \ell(u(\cdot, \bm y)) .
\end{equation*}
We are interested in evaluating the expectation and variance of $s$ as
\begin{equation*}
 E[s] = \int_{\Lambda} s( \bm y) \rho( \bm y) d \bm y,\qquad \qquad 
V[s] = \int_{\Lambda} \lp E[s] - s( \bm y)\rp^2 \rho( \bm y) d \bm y,
\end{equation*}
where $\rho( \bm y) = \prod_{q = 1}^Q \rho_q(y_q)$. Below we describe the hybridizable discontinuous Galerkin method for solving the model problem (\ref{eq:helm}) and the Monte Carlo simulation for computing estimates of $E[s]$ and $V[s]$.

\subsection{HDG discretization}

We begin by rewriting the governing equation (\ref{eq:helm}) as a first-order system
\begin{subequations}\label{eq:theproblem}
\begin{alignat}{2}
\boldsymbol{q} - \ggrad u & = 0,  \quad\mbox{ in }\mathcal{D}, \\
- \nabla\cdot\, \kappa \boldsymbol{q}   + \varrho \,u      & = f,
\quad\mbox{ in }\mathcal{D},\\
\kappa \boldsymbol{q}\cdot\boldsymbol{n}+ \nu u  & = g        \quad\mbox{ on }\partial \mathcal{D} .
\end{alignat}
\end{subequations}
The physical domain $\mathcal{D}$ is triangulated into elements $T$
forming a mesh $\mathcal{T}_h$ satisfying the standard finite element
conditions \cite{ciarlet1978finite}. Then,  letting  ${\partial{\mathcal{T}_h}}:=\{\partial T: T\in\mathcal{T}_h\}$
and denoting by  $\mathcal{F}_h$ the set of the faces $F$ of the elements $T\in\mathcal{T}_h$,
we seek  a vector approximation $\boldsymbol{q}_h \in \boldsymbol{V}_h^p$ to $\boldsymbol{q}$, a scalar approximation $u_h \in  W_h^p$ to $u$, 
and {a scalar approximation $\widehat{u}_h \in M_h^p$ to the {\em trace} of $u$  {on}  element {boundaries}}, where 
\begin{alignat*}{2}
  \boldsymbol{V}_h^p &=  \{ \boldsymbol{v}\in\boldsymbol{L}^2(\mathcal{D}):  \boldsymbol{v}|_T \in [P_p(T)]^d&&\;\forall T\in\mathcal{T}_h\},
  \\
  W_h^p &= \{ w\in L^2(\mathcal{D}): w|_T \in P_p(T) && \;\forall T\in\mathcal{T}_h\},
  \\
  M_h^p &= \{ \mu\in L^2(\mathcal{F}_h): \mu|_F \in
  P_p(F)&&\;\forall F\in\mathcal{F}_h\},
\end{alignat*}
and $P_p(D)$ is a space  of complex-valued polynomials of degree at most $p$ on $D$. Note that $\widehat{u}_h$ are defined only on the faces of the elements, hence they are single valued. We introduce the following inner products 
\begin{equation*}
(v,{w})_{{\mathcal{T}_h}}:=\sum_{T \in \mathcal{T}_h} (v,{w})_T, \qquad \langle{v,{w}}\rangle_{{\partial{\mathcal{T}_h}}} :=
\sum_{T \in \mathcal{T}_h} \langle{ v,{w}}\rangle_{\partial T},
\end{equation*}
where $(u,{v})_D :=\int_D u \overline{v} \; dx$ whenever $D$ is a domain in
$\mathbb{R}^d$, and $\langle{u,{v}}\rangle_D := \int_D u \overline{v} \; dx$ whenever $D$ is a domain
in $\mathbb{R}^{d-1}$.  For vector-valued functions $\boldsymbol{v}$ and $\boldsymbol{w}$, the integrals
are similarly defined with the integrand being the dot product $\boldsymbol{v}
\cdot \overline{\boldsymbol{w}}$. Note that  $\overline{w}$ denotes the complex conjugate of $w$.

The HDG approximations $(\boldsymbol{q}_h,u_h,\widehat{u}_h)$ in $\boldsymbol{V}_h^p \times W_h^p \times  M_h^p$
are determined by requiring that
\begin{subequations}\label{eq0:method}
\begin{alignat}{2}
 (\boldsymbol{q}_h,{\boldsymbol{r}})_{{\mathcal{T}_h}} +
 ( u_h,\nabla\cdot {\boldsymbol{r}})_{{\mathcal{T}_h}} -   
 \langle{\widehat{u}_h,{\boldsymbol{r}}\cdot\boldsymbol{n}}\rangle_{{\partial{\mathcal{T}_h}}}
 &  = 0 \label{eq0:method-a},
 \\
(\kappa \boldsymbol{q}_h, \ggrad {w})_{{\mathcal{T}_h}} -
 \langle{ \kappa \widehat{\boldsymbol{q}}_h\cdot\boldsymbol{n}, {w}}\rangle_{{\partial{\mathcal{T}_h}}}
 +(\varrho u_h,{w})_{\mathcal{T}_h}
 & =
 (f, {w})_{{\mathcal{T}_h}},  \label{eq0:method-b}
\\
\langle{\kappa \widehat{\boldsymbol{q}}_h \cdot \boldsymbol{n},{\mu} }\rangle_{\partial\mathcal{T}_h} + \langle{\nu \widehat{u}_h,{\mu} }\rangle_{\partial\mcal{D}}
& = \langle{g, {\mu}}\rangle_{\partial\mathcal{D}},  \label{eq0:method-c}
\end{alignat}
\end{subequations}
hold for all $(\boldsymbol{r},w,\mu)$ in $\boldsymbol{V}_h^p \times W_h^p \times M_h^p$, where the numerical flux ${\widehat{\boldsymbol{q}}}_h$ is defined as
\begin{equation}
\label{flux}
{\widehat{\boldsymbol{q}}}_h = \boldsymbol{q}_h - \tau \, \big( u_h - \widehat{u}_h \big) \boldsymbol{n}, \qquad\mbox{ on } \partial\mathcal{T}_h .
\end{equation}
Here  $\tau$ is the so-called {\em stabilization} parameter, a global constant with dimensions $\tau = 1/L$ where $L$ is the characteristic lengthscale. We set $\tau = 1$ since we do not consider multiple physical scales in this work. Further discussions on $\tau$ may be found in \cite{cockburn2009unified,nguyen2009linearCD}. By substituting (\ref{flux}) into (\ref{eq0:method}) we obtain that $(\boldsymbol{q}_h,u_h,\widehat{u}_h) \in \boldsymbol{V}_h^p \times W_h^p \times M_h^p$ satisfies
\begin{subequations}\label{method2}
\begin{alignat}{2}
 (\boldsymbol{q}_h,{\boldsymbol{r}})_{{\mathcal{T}_h}} +
 ( u_h,\nabla\cdot {\boldsymbol{r}})_{{\mathcal{T}_h}} -   
 \langle{\widehat{u}_h,{\boldsymbol{r}}\cdot\boldsymbol{n}}\rangle_{{\partial{\mathcal{T}_h}}}
 &  = 0,\label{eq0:method2-a}
 \\
(\kappa \boldsymbol{q}_h, \ggrad {w})_{{\mathcal{T}_h}} -
 \langle{ \kappa {\boldsymbol{q}}_h\cdot\boldsymbol{n} - \kappa \tau ( u_h - \widehat{u}_h ), {w}}\rangle_{{\partial{\mathcal{T}_h}}}
 +(\varrho u_h,{w})_{\mathcal{T}_h}
 & =
 (f, {w})_{{\mathcal{T}_h}},  \label{eq0:method2-b}
\\
 \langle{\kappa {\boldsymbol{q}}_h \cdot \boldsymbol{n}  - \kappa \tau ( u_h - \widehat{u}_h ),\mu\rangle_{\partial\mathcal{T}_h} + \langle\nu \widehat{u}_h,{\mu} }\rangle_{\partial\mathcal{D}}
& = \langle{g, {\mu}}\rangle_{\partial\mathcal{D}}, \label{eq0:method2-c}
\end{alignat}
\end{subequations}
for all $(\boldsymbol{r},w,\mu)$ in $\boldsymbol{V}_h^p \times W_h^p \times M_h^p$. This completes the definition of the HDG method.

The above weak formulation of the HDG method involves three field variables, namely, $\boldsymbol{q}_h,u_h,$ and $\widehat{u}_h$.  However, the first two equations (\ref{eq0:method2-a}) and (\ref{eq0:method2-b}) allow us to write both $\bm{q}_h$ and $u_h$ in terms of $\widehat{u}_h$ at the element level due to the fact that our approximation spaces are discontinuous. Therefore, we can substitute both $\bm{q}_h$ and $u_h$ from the first two equations into the last equation (\ref{eq0:method2-c}) to obtain a weak formulation in terms of $\widehat{u}_h$ only: find $\widehat{u}_h \in M_h^p$ such that
\begin{equation}
\label{reducedform}
\widehat{a}_h(\widehat{u}_h, \mu) = \widehat{b}_h(\mu), \qquad \forall \, \mu \, \in M_h^p .
\end{equation} 
Here we omit the derivation of the bilinear form $\widehat{a}_h$ and the linear functional $\widehat{b}_h$. Instead we refer the reader to \cite{nguyen2009linearCD} for a detailed discussion. The reduced weak formulation (\ref{reducedform}) gives rise to the following linear system
\begin{equation}
\label{matrixsystem}
\widehat{\mathbf{A}} \widehat{\mathbf{u}} = \widehat{\mathbf{b}} ,
\end{equation}
where $\widehat{\mathbf{u}}$ is the vector containing the degrees of freedom of $\widehat{u}_h$. Because $\widehat{u}_h$ is single valued on the faces of the finite element mesh, it has significantly fewer degrees of freedom than $u_h$. As a result, the global matrix system (\ref{matrixsystem}) of the HDG method can be much smaller than that of other DG methods. This results in significant savings in terms of computational time and memory storage. 

 It turns out that although the formulation (\ref{reducedform}) results in the smallest possible system, it is not ideal to use it as the starting point for our reduced basis method. Substituting the first two equations (\ref{eq0:method2-a}) and (\ref{eq0:method2-b}) into the last equation (\ref{eq0:method2-c}) results in the inverse of the material coefficients $\kappa$ and $\varrho$, which renders the bilinear form $\widehat{a}_h$  {\em nonaffine} in the material coefficients. Although nonaffine parameter dependence can be treated by using the empirical interpolation method \cite{barrault2004empirical} or the best points 
interpolation method \cite{nguyen2008best}, such treatment incurs additional cost and is unnecessary. We are going to derive a new weak formulation of the HDG method, which is suited for the reduced basis method.

\subsection{A new weak formulation of the HDG method}

We begin by deriving a weak formulation of the HDG method upon which our reduced basis method is constructed. To this end, we introduce two lifting operators $\bm l : W_h^p \to \bm V_h^p$ and $\bm m : M_h^p \to \bm V_h^p$ defined as
\begin{subequations}
\label{qeq}
\begin{alignat}{2}
 (\bm l (w),{\boldsymbol{r}})_{{\mathcal{T}_h}}  & =  -( w,\nabla\cdot {\boldsymbol{r}})_{{\mathcal{T}_h}}, &\qquad  \forall \, \bm r \in \bm V_h^p, \\
 (\boldsymbol{m}(\mu),{\boldsymbol{r}})_{{\mathcal{T}_h}} & =   \langle{\mu,{\boldsymbol{r}}\cdot\boldsymbol{n}}\rangle_{{\partial{\mathcal{T}_h}}}, & \qquad  \forall \, \bm r \in \bm V_h^p .
\end{alignat}
\end{subequations}
It thus follows from (\ref{eq0:method2-a}) and (\ref{qeq}) that we can express $\bm{q}_h$ as a function of $u_h$ and $\widehat{u}_h$ as
\begin{equation} 
\label{qlift}
\bm{q}_h = \bm l (u_h) + \boldsymbol{m}({\widehat{u}}_h) .
\end{equation}
By substituting (\ref{qlift}) into (\ref{eq0:method2-b}) and (\ref{eq0:method2-c}) we arrive at the following weak formulation: find $(u_h,\widehat{u}_h) \in W_h^p \times M_h^p$ such that
\begin{align*}
\bigl( \kappa (\bm l (u_h) + \boldsymbol{m}({\widehat{u}}_h)), \ggrad {w}\bigr)_{{\mathcal{T}_h}}  -
 \bigl\langle{ \kappa (\bm l (u_h) + \boldsymbol{m}({\widehat{u}}_h)) \cdot\boldsymbol{n} - \kappa \tau ( u_h - \widehat{u}_h ), {w}}\bigr\rangle_{{\partial{\mathcal{T}_h}}}+(\varrho u_h,{w})_{\mathcal{T}_h} &=
 (f, {w})_{{\mathcal{T}_h}},
\\
\bigl\langle{\kappa (\bm l (u_h) + \boldsymbol{m}({\widehat{u}}_h)) \cdot \boldsymbol{n}  - \kappa \tau ( u_h - \widehat{u}_h ),{\mu} }\bigr\rangle_{\partial\mathcal{T}_h} + \langle \nu \widehat{u}_h,\mu \rangle_{\partial\mathcal{D}} &= \langle{g, {\mu}}\rangle_{\partial\mathcal{D}},
\end{align*}
for all $(w,\mu)$ in $W_h^p \times M_h^p$. 
By setting the $\mcal{N}$-dimensional approximation space to be $\bm W_h^p := W_h^p \times M_h^p$, $\bm u_h := (u_h,\widehat{u}_h)$, and $\bm{w} := (w,\mu)$ we obtain that $\bm u_h \in \bm W_h^p$, satisfies
\begin{equation}
\label{weakuuh}
a_h(\bm{u}_h, \bm{w}; (\kappa,\varrho,\nu)) = b_h(\bm{w}), \qquad \forall \, \bm{w} \in \bm W_h^p,
\end{equation}
where the bilinear form $a_h$ and the linear functional $b_h$ are given by
\begin{subequations}
\label{aform}
\begin{alignat}{2}
a_h(\bm{v}, \bm{w}; (\kappa,\varrho,\nu))  = &  \bigl(\kappa \lp\bm l(v) + \boldsymbol{m}(\eta)\rp, \ggrad {w}\bigr)_{{\mathcal{T}_h}} - \bigl\langle{ \kappa (\bm l(v) + \boldsymbol{m}(\eta)) \cdot\boldsymbol{n} - \kappa \tau ( v - \eta ), {w}}\bigr\rangle_{{\partial{\mathcal{T}_h}}} \nonumber \label{aforma}\\& +(\varrho v,{w})_{\mathcal{T}_h}  + \,  \bigl\langle{\kappa (\bm l(v) + \boldsymbol{m}(\eta)) \cdot \boldsymbol{n}  - \kappa \tau(v - \eta),{\mu} }\bigr\rangle_{\partial\mathcal{T}_h} + \langle \nu \eta,\mu \rangle_{\partial \mcal{D}}, \\[1ex]
b_h(\bm{w}) =  & \ (f, {w})_{{\mathcal{T}_h}} + \langle{g, {\mu}}\rangle_{\partial\mathcal{D}}, 
\end{alignat}
\end{subequations}
for all $\bm{v} := (v,\eta) \in \bm W_h^p$ and $\bm{w} := (w,\mu) \in \bm W_h^p$. We note that the bilinear form \eqref{aforma} is affine in $\bm y = (\kappa,\varrho,\nu)$. Furthermore, if we select $\bm r = \kappa \lp\bm l(v) + \bm m(\eta)\rp$ in \eqref{qeq} and substitute into \eqref{aforma}, we recover a symmetric form, which is also coercive provided $\kappa,\varrho,\nu>0$. Henceforth, the choice $(\kappa,\varrho,\nu)= (1,1,1)$ allows us to equip the approximation space $\bm W_h^p$ with the inner product $(\bm{v}, \bm{w})_{\bm W} := a_h(\bm{v}, \bm{w}; (\kappa,\varrho,\nu)= (1,1,1))$ and the induced norm  $\|\bm w\|_{\bm W} = \sqrt{(\bm{w}, \bm{w})_{\bm W}}$.

We now substitute the expression of $\kappa$ from (\ref{affinekappa}) into (\ref{aform}) to express $a_h$ as
\begin{equation}
\label{affineform}
a_h\lp\bm{v}, \bm{w}; (\kappa,\varrho,\nu)\rp = a_h(\bm{v}, \bm{w}; \bm y) = a_h^0(\bm{v}, \bm{w}) + \sum_{q=1}^Q y_q a_h^q(\bm{v}, \bm{w}) ,
\end{equation}
where the bilinear forms are given by $a_h^q(\bm{v}, \bm{w}) := a_h\lp\bm{v}, \bm{w}; (\psi_q,0,0)\rp$ for $1 \le q \le Q $ and $a_h^0(\bm{v}, \bm{w}) := a_h\lp\bm{v}, \bm{w}; (\overline{\kappa},\varrho,\nu)\rp$. Therefore, we can write the weak formulation (\ref{weakuuh}) as follows: for any $\bm y \in \Lambda$, $\bm u_h(\bm y) \in \bm W_h^p$ satisfies
\begin{equation}
\label{weakform}
a_h(\bm{u}_h, \bm{w}; \bm y) = b_h(\bm{w}), \qquad \forall \, \bm{w} \in \bm W_h^p .
\end{equation}
Finally, we evaluate our realization output as
\begin{equation*}
s_h(\bm y) = \ell_h(\bm{u}_h(\bm y)) ,
\end{equation*}
where the linear functional $\ell_h$ is obtained from the HDG discretization of $\ell$. The key point of the new HDG formulation \eqref{weakform} for an efficient perfomance of the reduced basis method is the affine representation \eqref{affineform}. This aspect is of crucial importance, and the main reason we prefer \eqref{weakform} to the reduced weak formulation \eqref{reducedform} for constructing the reduced basis approximation. Furthermore, the new formulation is optimal in terms of degrees of freedom, since we no longer account for the gradient $\bm{q}_h$. Finally, even though the parameter independent matrices arising from \eqref{affineform} are used for the reduced basis approximation, the solution $\bm u_h$ is never computed as the solution of the full system \eqref{weakform}. Instead, we can invoke again discontinuity of the approximation spaces to write $u_h$ in terms of $\widehat{u}_h$. This common strategy in HDG methods enables us to solve for the global degrees of freedom of $\widehat{u}_h$ only and then recover $u_h$ efficiently.
 \subsection{Monte Carlo sampling with the HDG method}
We are interested in evaluating statistics of the output $s_h(\bm y)$ such as its expectation and variance. Let $Y_M = \{\bm{y}_m  \in \Lambda, 1 \le m \le M\}$ be a set of random samples drawn in the parameter space $\Lambda$ with the probability density function $\rho(\bm y)$. We evaluate the following outputs 
\begin{equation}\label{hdgoutputs}
s_h(\bm y_m) = \ell_h(\bm u_h(\bm y_m)), \qquad m = 1,\ldots, M .
\end{equation}
The Monte Carlo-HDG (MC-HDG) estimates of the expectation $E[s]$ and variance $V[s]$ can be computed, respectively, as
\begin{equation}
\label{hdgmc}
E_M[s_h] =  \frac{1}{M} \sum_{m=1}^M s_h(\bm y_m), \qquad V_M[s_h] =  \frac{1}{M-1} \sum_{m=1}^M \lp E_M[s_h] - s_h(\bm y_m)\rp^2 .
\end{equation}
We shall assume that $s_h(\bm y)$ is indistinguishable from $s(\bm y)$ for any $\bm y \in \Lambda$. Moreover, it is a known result that the estimators in \eqref{hdgmc} are unbiased and converge in distribution to
\beq
E[s_h] - E_M[s_h] \xrightarrow{\; d} N\lp 0\,; \frac{V[s_h]}{M} \rp,\qquad V[s_h] - V_M[s_h] \xrightarrow{\; d} N\lp 0\,; \frac{V[\lp s_h - E[s_h]\rp^2]}{M} \rp\:.
\eeq
Confidence intervals can be constructed employing the central limit theorem (CLT), that is for all $a>0$ we have
\begin{subequations}
\begin{align}
\lim_{M \to \infty} \mathrm{Pr}\left(\abs{E[s_h] - E_M[s_h]} \le \Delta_{h,M}^E \right)  &=  \mathrm{erf} \left(\frac{a}{\sqrt{2}} \right) ,\label{MCHDGerror}\\
\lim_{M \to \infty} \mathrm{Pr}\left(\abs{V[s_h] - V_M[s_h]} \le \Delta_{h,M}^V \right)  &=  \mathrm{erf} \left(\frac{a}{\sqrt{2}} \right) ,
\end{align}
\end{subequations}
where
\begin{equation}
\label{MCHDGbound}
\Delta_{h,M}^E = a \sqrt{\frac{V_M[s_h]}{M}} \:,\qquad \Delta_{h,M}^V = a \sqrt{\frac{V_M[\lp s_h - E_M[s_h]\rp^2]}{M}}\:.
\end{equation}
Therefore, in order to guarantee that $\abs{E[s_h] - E_M[s_h]}$ is bounded by a specified error tolerance $\epsilon_{\rm tol}$ with a high probability (say, greater than $0.95$), we need to take $a \ge 1.96$ and $M \ge a^2 V_M[s_h]/ \epsilon_{\rm tol}^2$.  As a result, $M$ can be very large for a small error tolerance. Hence, the evaluations \eqref{hdgoutputs}--\eqref{hdgmc} can be very demanding.

The remaining goals of this paper are as follows. On one hand, we develop a reduced basis (RB) method for rapid reliable approximation of the stochastic HDG output $s_h(\bm y)$ for any given parameter vector $\bm{y}$ in $\Lambda$. On the other hand, we develop a multilevel variance reduction method to accelerate the convergence of the Monte Carlo simulation by exploiting the exponentially fast convergence of the RB output to the high-fidelity HDG output as a function of the RB dimension. These two ingredients enable very fast reliable computation of the statistical outputs at a computational cost which is several orders of magnitude less expensive than that of the MC-HDG approach. We describe the reduced basis approach in Section 3 and the multilevel variance reduction method in Section 4.

\section{Reduced Basis Method}\label{sec:RB}
We consider a ``primal-dual'' formulation \cite{pierce2000adjoint,cuong2005certified} particularly well-suited to good approximation and error characterization of the output. To this end, we introduce the dual problem of \eqref{weakform}:  given $\bm y \in \Lambda$, the dual solution $\bm \phi_h(\bm y) \in \bm W_h^p$ satisfies
\begin{equation*}
  a_h(\bm v,\bm \phi_h;\bm y) = - \ell_h(\bm v), \qquad \forall \ \bm v \in \bm W_h^p \ .
\end{equation*}
The dual problem plays an important role in improving the convergence rate of both the RB output and associated error bound. 

We next assume that we are given orthonormalized basis functions $\bm \zeta_n^{\rm pr}, \bm \zeta_n^{\rm du}  \in \bm W_h^p$, $1 \le n \le N_{\max}$, such that $(\bm \zeta_m^{\rm pr}, \bm \zeta_n^{\rm pr})_{\bm W} = (\bm \zeta_m^{\rm du}, \bm \zeta_n^{\rm du})_{\bm W} = \delta_{mn}, 1 \le m,n \le N_{\max}$. We define the associated hierarchical RB spaces as
\begin{equation*}
\bm W^{\rm pr}_N = \mbox{span} \{\bm \zeta_n^{\rm pr}, 1 \le n \le N\}, \quad \bm W^{\rm du}_N = \mbox{span} \{\bm \zeta_n^{\rm du}, 1 \le n \le N\}, \quad N=1, \ldots, N_{\max} \ .
\end{equation*}
In practice, the spaces $\bm W_N^{\rm pr}$ and $\bm W_N^{\rm du}$ consist of orthonormalized primal and dual solutions $\bm \zeta_n^{\rm pr}, \bm \zeta_n^{\rm du}$ at selected parameter values generated by a Greedy sampling procedure \cite{rozza2008reduced,grepl2005reduced,veroy2003posteriori}. For our present purpose, however, $\bm W_N^{\rm pr}$ and $\bm W_N^{\rm du}$ can in fact represent any sequence of (low-dimensional) hierarchical approximation spaces.  We then apply the Galerkin projection for both the primal and dual problems: Given $\bm y \in \Lambda$, we find a primal RB approximation $\bm u_{N}(\bm y) \in \bm W_N^{\rm pr}$ satisfying
\begin{equation}
\label{eq:RBprimal}
a_h(\bm u_{N}(\bm y),\bm w;\bm y) = b_h(\bm w), \qquad \forall\: \bm w \in \bm W_N^{\rm pr} \: ,
\end{equation}
and a dual RB approximation $\bm \phi_{N}(\bm y) \in \bm W_N^{\rm du}$ satisfying
\begin{equation*}
a_h(\bm w, \bm \phi_{N}(\bm y);\bm y) = -\ell_h(\bm w), \qquad \forall\: \bm w \in \bm W_N^{\rm du} \: .
\end{equation*}
We can now evaluate the RB realization output as
\begin{equation*}
s_{N}(\bm y) =  \ell_h(\bm u_N(\bm y)) + a_h(\bm u_N(\bm y),\bm \phi_N(\bm y); \bm y) - b_h(\bm \phi_N(\bm y)) .
\end{equation*}
As discussed below, the online computational cost of evaluating the RB output depends only on $N$ and $Q$. Hence, for small $N$ and $Q$, the RB approximation can be significantly less expensive than the HDG approximation.

The RB output is then used as an approximation to the HDG output in the Monte Carlo simulation.  The Monte Carlo-Reduced Basis (MC-RB) estimates of the expectation and variance of the output of interest are given by
\begin{equation*}
E_M[s_{N}]  =  \frac{1}M \sum_{m=1}^M s_{N}(\bm y_m), \qquad V_M[s_{N}]  =  \frac{1}{M-1} \sum_{m=1}^M \lp E_M [s_{N}] - s_{N}(\bm y_m) \rp^2 \,
\end{equation*}
for the same set of samples $Y_M = \{\bm{y}_m  \in \Lambda, 1 \le m \le M\}$. Since the RB approximation is constructed upon the HDG approximation these quantities actually approximate the MC-HDG estimates. We next develop {\em a posteriori} error bounds for our MC-RB estimates relative to the MC-HDG estimates.  

\subsection{\emph{A posteriori} error estimation}

We note from~\eqref{eq:RBprimal} that the residuals $r_h^{\rm pr}(\bm w;\bm y)$ and $r_h^{\rm du}(\bm w;\bm y)$ associated with $\bm u_{N}(\bm y)$ and $\bm \phi_{N}(\bm y)$, respectively, are given by
\begin{equation*}
r_h^{\rm pr}(\bm w;\bm y) = b_h(\bm w) - a_h(\bm u_{N}(\bm y),\bm w;\bm y), \quad r_h^{\rm du}(\bm w;\bm y) = -\ell_h(\bm w) - a_h(\bm w, \bm \phi_{N}(\bm y);\bm y),
\end{equation*}
for all $\bm w \in \bm{W}_h^p$.  The dual norm of the primal residual and the dual norm of the dual residual are given by
\begin{equation*}
\|r_h^{\rm pr}(\cdot;\bm y)\|_{\bm{W}'} = \sup_{\bm w \in \bm{W}_h^p} \frac{r_h^{\rm pr}(\bm w;\bm y)}{\|\bm w\|_{\bm W}}, \qquad \|r_h^{\rm du}(\cdot;\bm y)\|_{\bm{W}'} = \sup_{\bm w \in \bm{W}_h^p} \frac{r_h^{\rm du}(\bm w;\bm y)}{\|\bm w\|_{\bm W}} \ .
\end{equation*}
It is a standard result~\cite{cuong2005certified,rozza2008reduced} that
\begin{alignat*}{1}
\|\bm u_h(\bm y) - \bm u_{N}(\bm y)\|_{\bm W} & \le  \Delta^{\rm pr}_{N}(\bm y) \equiv \frac{\|r_h^{\rm pr}(\: \cdot \:;\bm y)\|_{\bm W'}}{\tilde{\beta}(\bm y)} \ ,
\\
\|\bm \phi_h(\bm y) - \bm \phi_{N}(\bm y)\|_{\bm W} & \le  \Delta^{\rm du}_{N}(\bm y) \equiv \frac{\|r_h^{\rm du}(\: \cdot \:;\bm y)\|_{\bm W'}}{\tilde{\beta}(\bm y)} \ , \\
 |s_h(\bm y)-s_N(\bm y)| & \le \Delta_N^s(\bm y) \equiv \tilde{\beta}(\bm y) \Delta_N^{\rm pr}(\bm y) \Delta_N^{\rm du}(\bm y) \ ,
\end{alignat*}
where $\tilde{\beta}(\bm y)$ is a positive lower bound for the Babu\v{s}ka ``inf-sup'' stability constant $\beta_h(\bm y)$ defined as
\begin{equation*}
0<\beta_h(\bm y) \equiv \inf_{\bm w \in \bm W_h^p} \: \sup_{\bm v \in \bm W_h^p}\:  \frac{a_h(\bm w, \bm v;\bm y)}{ \| \bm w\|_{\bm W} \| \bm v \|_{\bm W}},
\end{equation*}
that is, the minimum (generalized) singular value associated with the differential operator. It is critical to note that  the output error (and output error bound) vanishes as the {\em product\/} of the primal and dual error (bounds), and hence much more rapidly than either the primal or dual error.  

It thus follows that we can bound the errors in the MC-RB estimates relative to the MC-HDG estimates as
\begin{equation}
 |E_M[s_h] - E_M[s_{N}] | \le \frac{1}M \sum_{m=1}^M |s_{h}(\bm y_m) - s_{N}(\bm y_m)| \le  \frac{1}M \sum_{m=1}^M \Delta_N^s(\bm y_m) \equiv \Delta_{N,M}^{E},
   \label{MCRBExp}
\end{equation}
and
{\begin{eqnarray}
 \abs{V_M[s_h] - V_M[s_{N}]} & = &  \frac{1}{M-1} \left\lvert\sum_{m=1}^M \Bigl( E_M [s_{h}] - s_{h}(\bm y_m) \Bigr)^2 - \Bigl( E_M [s_{N}] - s_{N}(\bm y_m) \Bigr)^2\right\rvert \nonumber\\
 & = & \frac{1}{M-1} \Bigg\lvert\sum_{m=1}^M  \Bigl( s_{h}(\bm y_m)  - s_{N}(\bm y_m) - E_M [s_{h}] + E_M [s_{N}] \Bigr)  \Bigl( s_{h}(\bm y_m)   + s_{N}(\bm y_m) \Bigr)\Bigg\rvert \nonumber \\
 & \le & \frac{1}{M-1} \sum_{m=1}^M  \Bigl(\abs{s_{h}(\bm y_m)  - s_{N}(\bm y_m)} + \abs{E_M [s_{h}] - E_M [s_{N}]} \Bigr) \abs{s_{h}(\bm y_m)   + s_{N}(\bm y_m)} \nonumber \\
  & \le & \frac{1}{M-1} \sum_{m=1}^M \Bigl(\Delta_N^s(\bm y_m) + \Delta_{N,M}^{E}\Bigr) \Bigl(\Delta_N^s(\bm y_m)  + 2\abs{s_{N}(\bm y_m)} \Bigr) \equiv \Delta_{N,M}^{V} .
  \label{MCRBVar}
\end{eqnarray}}

It should be stated that this error bound is rather pessimistic, and that a more precise bound can be obtained by introducing suitable dual problems to recover a quadratically convergent bound for the variance, as reported in \cite{haasdonk2013reduced}. We can also bound the difference between the RB expected value and the true value. To this end, we note from the triangle inequality that
\begin{equation}\label{MCRBerror}
\abs{E[s_h] - E_M[s_{N}]}  \le  \abs{E[s_h] - E_M[s_{h}] }+ \abs{E_M[s_h] - E_M[s_{N}]} ,
\end{equation}
Following from (\ref{MCHDGerror}), (\ref{MCHDGbound}), (\ref{MCRBExp}), (\ref{MCRBVar}), and (\ref{MCRBerror}) we define the error bound
\begin{equation}\label{eq:CLTMCRBbound}
\widetilde{\Delta}_{N,M}^{E} = a \sqrt{\frac{(V_M[s_N] + \Delta_{N,M}^{V} )}{M}} + \Delta_{N,M}^{E}\\
\end{equation}
such that 
\begin{equation*}
\lim_{M \to \infty} {\rm Pr} \left (\abs{E[s_h] - E_M[s_{N}]} \le \widetilde{\Delta}_{N,M}^{E} \right) \geq \mathrm{erf}\left(\frac{a}{\sqrt{2}} \right).
\end{equation*}
Clearly, the error bound (\ref{eq:CLTMCRBbound}) comprises two terms: the first term is due to the MC sampling, while the second term is due to the RB approximation.

\subsection{Computational strategy}

The linearity and parametric affinity of the problem allow for an efficient Offline-Online decomposition strategy. The Offline stage --- parameter independent, computationally intensive but performed only once --- comprises the greedy search for the selection of parameter values, the computation of snapshots $\bm \zeta_n^{\rm pr}, \bm \zeta_n^{\rm du}, 1 \le n \le N_{\max}$ associated with the HDG approximation space at the selected parameter values and the formation and storage of several parameter-independent small matrices and vectors. The Online stage --- parameter dependent, performed multiple times --- evaluates $s_{N}(\bm y),\,\Delta^s_N(\bm y)$ for any new $\bm y$ with complexity $\mcal{O}\lp 2N^3 + 2(Q+1)^2N^2\rp$ independent of the dimension $\mcal{N}$ of the HDG approximation space. The implications are twofold: first, if $N$ and $Q$ are indeed small, we shall achieve very fast output evaluation, usually several orders of magnitude faster than the HDG output; second, we may choose the 
HDG 
approximation very conservatively --- to effectively eliminate the error between the exact output and HDG  output --- without adversely affecting the Online (marginal) cost. We refer the reader to \cite{cuong2005certified,prud2002reliable} for a more thorough description of the Offline-Online procedure.

It is clear that the error in the RB expected value and its error bound depend on $N$ and $Q$ as well as on $M$. Typically, both the error and its error bound decrease very rapidly as a function of $N$, but very slowly as a function of the number of samples $M$. Hence, $M$ should be chosen very large, while $N$ can be chosen to be much smaller. Indeed, the (Online) computational cost to evaluate the RB expected value $E_M[s_{N}]$ and its error bound $\widetilde{\Delta}_{N,M}^{E} $ scales as $\mcal{O}\lp2M(N^3+2(Q+1)^2N^2)\rp$.  Since both $Q$ and $N$ are typically very small, the RB method can effect significant savings relative to the HDG method. Nevertheless, its performance can be affected by the accuracy of the RB outputs and the sharpness of the RB error bounds.

\section{Model and Variance Reduction Method}\label{sec:RBCV}

\subsection{Control variates principle}\label{ssec:CV}

We first review the essential idea of control variates, which will serve as a building block for our method. Let $X$ be a random variable. We would like to estimate the expected value of $X$. Suppose that we have another  random variable $Y$ and that its expected value $E[Y]$ is either known or inexpensive to compute. We then introduce a new random variable
\begin{equation*}
X^* = X + \gamma (E[Y] - Y) ,
\end{equation*}
where $\gamma$ is a deterministic coefficient. It is obvious that $E[X^*] = E[X]$ for any choice of $\gamma$. However, the variance of $X^*$ is different from that of $X$. Specifically, we have
\begin{equation*}
V[X^*] = V[X] + \gamma^2 V[Y] - 2 \gamma \mathrm{Cov}(X,Y) ,
\end{equation*}
where $ \mathrm{Cov}(X,Y) = E[XY] - E[X]E[Y]$ is the covariance of $X$ and $Y$. It can be easily shown that the following choice 
\begin{equation*}
\gamma = \frac{ \mathrm{Cov}(X,Y)}{V[Y]}
\end{equation*}
is optimal in the sense that it minimizes the variance of $X^*$. With this choice, we have
\begin{equation*}
V[X^*] = V[X] (1 - \rho^2(X,Y)), 
\end{equation*}
where $\rho(X,Y) =  \mathrm{Cov}(X,Y)/\sqrt{V[X] V[Y]}$ is the correlation coefficient of $X$ and $Y$. It is clear that if $X$ and $Y$ are highly correlated (i.e., $\rho(X,Y)$ is close to $\pm 1$) then $V[X^*]$ is much smaller than $V[X]$. In that case, the MC simulation of $E[X^*]$ converges significantly faster than that of $E[X]$ according to the CLT. 

In summary, control variate methods try to estimate $E[X]$ by using the ``surrogate'' expected value $E[Y]$ and sampling the reduced variance variable $X^*$. When the same principle is applied recursively to estimate $E[Y]$, the resulting method is called multilevel control variates. 

\subsection{Two-level Monte Carlo sampling}

We now apply the above idea to compute an estimate of $E[s_h]$, where $s_h(\bm y)$ is the stochastic output obtained by using the HDG method to solve the underlying stochastic PDE as described in Section 2. To achieve this goal, we introduce 
\begin{equation*}
s_h^*(\bm y) = s_h(\bm y) + \gamma (E[s_{N_1}] - s_{N_1}(\bm y)),
\end{equation*}
where $s_{N_1}(\bm y)$ is the RB output developed in Section 3 for some $N_1 \in [1, N_{\max}]$. Because $s_{N_1}(\bm y)$ generally approximates $s_h(\bm y)$ very well, the two outputs are highly correlated. Therefore, we choose $\gamma = 1$ as we expect that the optimal value of $\gamma$ is close to 1. With this choice, we obtain
\begin{equation}\label{eq:2LE_exact}
E[s_h] = E[s_h^*] =  E[s_h - s_{N_1}] + E[s_{N_1}] .
\end{equation}
The underlying premise here is that the two expectation terms on the right hand side can be computed efficiently by MC simulations owing to variance reduction and model reduction: the first term requires a small number of samples because its variance is generally very small, while the second term is less expensive to evaluate because it involves the RB output.

In particular, let $Y^0_{M_0} = \{\bm y^0_m \in \Lambda,  1 \le m \le M_0\}$ and $Y^1_{M_1} = \{\bm y^1_m \in \Lambda,  1 \le m \le M_1\}$ be two independent sets of random samples drawn in $\Lambda$ with the probability density function $\rho(\bm y)$. We calculate our Model and Variance Reduction (MVR) unbiased estimate of $E[s_h]$ as
\begin{equation}\label{eq:2LE}
E_{M_0,M_1}[s_h] =  E_{M_0}[s_h -s_{N_1}]  + E_{M_1}[s_{N_1}] ,
\end{equation}
where
\begin{equation}\label{eq:2LE_MC}
E_{M_0}[s_h - s_{N_1}]  = \frac{1}{M_0} \sum_{m=1}^{M_0} \lp s_h(\bm y^0_m) - s_{N_1}(\bm y^0_m)\rp, \qquad E_{M_1}[s_{N_1}]   =  \frac{1}{M_1} \sum_{m=1}^{M_1} s_{N_1}(\bm y_m^1)
\end{equation}
We note that our approach computes an estimate of $E[s_h]$, while the MC-RB approach described in the previous section computes an estimate of $E[s_N]$.

Similarly, we exploit the control variates idea to compute an estimate of the true variance $V[s_h]$ given by
\begin{equation}\label{eq:2LV}
V_{M_0,M_1}[s_h] =  E_{M_0}[\zeta_h -\zeta_{N_1}]  + E_{M_1}[\zeta_{N_1}] \:,
\end{equation}
where $\zeta_h := \lp s_h - E_{M_0,M_1}[s_h] \rp^2$ and $\zeta_{N_1} := \lp s_{N_1} - E_{M_0,M_1}[s_h] \rp^2$ and the expectations in \eqref{eq:2LV} are analogous to the expectations in \eqref{eq:2LE_MC}. The variance estimate is negatively biased
\begin{equation*}
E\lb V_{M_0,M_1}[s_h] - V[s_h] \rb = -\frac{V[s_h - s_{N_1}]}{M_0} - \frac{V[s_{N_1}]}{M_1}\:,
\end{equation*}
as shown in the \ref{sec:bias}.

It remains to provide {\em a posteriori} estimates for the errors in the expectation and variance. Subtracting \eqref{eq:2LE} from \eqref{eq:2LE_exact} we identify new random variables $Z_0,\,Z_1$ whose limiting distributions are normal, and since they are independent their sum is also normally distributed,
\begin{subequations}\label{eq:RV_exp}
\begin{align}
Z_0 & = E[s_h - s_{N_1}] - E_{M_0}[s_h-s_{N_1}] \sim N\lp 0\,;\, \frac{V[s_h -s_{N_1}]}{M_0} \rp,\\
Z_1 & = E[s_{N_1}] - E_{M_1}[s_{N_1}] \sim N\lp 0\,;\, \frac{V[s_{N_1}]}{M_1} \rp,\\
Z_0 + Z_1 & = E[s_h] - E_{M_0,M_1}[s_h] \sim N\lp 0\,;\, \frac{V[s_h -s_{N_1}]}{M_0}+ \frac{V[s_{N_1}]}{M_1} \rp \label{eq:z1z2}.
\end{align}
\end{subequations}
We invoke now the CLT to obtain an error estimate for the expectation error as
\begin{equation}\label{eq:CLTE}
\lim_{M_0 \to \infty} \lim_{M_1 \to \infty} \mathrm{Pr} \left( \abs{E[s_h] - E_{M_0,M_1}[s_h]} \le \Delta_{M_0,M_1}^E \right) = \mathrm{erf}\left(\frac{a}{\sqrt{2}} \right)   ,
\end{equation}
where 
\begin{subequations}
\begin{align}
\Delta_{M_0,M_1}^E  & =  a \sqrt{\frac{ V_{M_0}[s_h - s_{N_1}]}{M_0}+ \frac{V_{M_1}[s_{N_1}]}{M_1}}, \label{eq:mlebound}\\
V_{M_0}[s_h - s_{N_1}] & = \frac{1}{M_0 - 1} \sum_{m=1}^{M_0} \lp E_{M_0}[s_h-s_{N_1}] - s_h(\bm y^0_m) + s_{N_1}(\bm y^0_m)\rp^2, \label{eq:mlev0}\\
V_{M_1}[s_{N_1}]  & =  \frac{1}{M_1-1} \sum_{m=1}^{M_1} \lp E_{M_1} [s_{N_1}] - s_{N_1}(\bm y_m^1)\rp^2 ,\label{eq:mlev1}
\end{align}
\end{subequations}
and the variances in \eqref{eq:z1z2} are estimated with their MC counterparts \eqref{eq:mlev0}--\eqref{eq:mlev1}.

For the variance, we first define auxiliary variables $\widehat{\zeta}_h := \lp s_h - E[s_h] \rp^2$ and $\widehat{\zeta}_{N_1} := \lp s_{N_1} - E[s_h] \rp^2$ and the auxiliary variance 
\beq
\widehat{V}_{M_0,M_1}[s_h] =  E_{M_0}[\widehat{\zeta}_h -\widehat{\zeta}_{N_1}]  + E_{M_1}[\widehat{\zeta}_{N_1}] \:.
\eeq
The MVR variance estimate in \eqref{eq:2LV} can be rewritten as (see \ref{sec:bias})
\beq
V_{M_0,M_1}[s_h] = \widehat{V}_{M_0,M_1}[s_h] - \bigl( E[s_h] - E_{M_0,M_1}[s_h] \bigr)^2,
\eeq
which implies
\beql\label{eq:VarCLt}
V_{M_0,M_1}[s_h] -V[s_h] = \bigl(\widehat{V}_{M_0,M_1}[s_h]-V[s_h]\bigr) - \bigl( E[s_h] - E_{M_0,M_1}[s_h] \bigr)^2 .
\eeql
Let us consider the two terms in the RHS in reverse order. Convergence in probability for the second term is guaranteed by \eqref{eq:CLTE}, that is $E_{M_0,M_1}[s_h] - E[s_h] \xrightarrow{\,P\,} 0$. Furthermore, repeating the analysis in \eqref{eq:RV_exp} for the first term leads to
\beql\label{eq:RV_var}
\widehat{V}_{M_0,M_1}[s_h]-V[s_h] \sim N\lp 0\,;\,\frac{V[\widehat{\zeta}_h-\widehat{\zeta}_{N_1}]}{M_0}+ \frac{V[\widehat{\zeta}_{N_1}]}{M_1}  \rp.
\eeql
Therefore the limiting distribution of $V_{M_0,M_1}[s_h] -V[s_h]$ is the same as the limiting distribution of $\widehat{V}_{M_0,M_1}[s_h]-V[s_h]$ (Slutzky's theorem), and the straightforward application of the CLT recovers
\begin{equation*}
\lim_{M_0 \to \infty} \lim_{M_1 \to \infty} \mathrm{Pr} \left( \abs{V[s_h] - V_{M_0,M_1}[s_h]} \le \Delta_{M_0,M_1}^V \right) = \mathrm{erf}\left(\frac{a}{\sqrt{2}} \right)   ,
\end{equation*}
where
\begin{subequations}
\begin{align}
\Delta_{M_0,M_1}^V  & =  a \sqrt{\frac{ V_{M_0}[\zeta_h - \zeta_{N_1}]}{M_0}+ \frac{V_{M_1}[\zeta_{N_1}]}{M_1}} \:, \label{eq:mlvbound}\\
V_{M_0}[\zeta_h - \zeta_{N_1}] & = \frac{1}{M_0 - 1} \sum_{m=1}^{M_0} \lp E_{M_0}[\zeta_h-\zeta_{N_1}] - \zeta_h(\bm y^0_m) + \zeta_{N_1}(\bm y^0_m)\rp^2, \label{eq:mlvv0}\\
V_{M_1}[\zeta_{N_1}]  & =  \frac{1}{M_1-1} \sum_{m=1}^{M_1} \lp E_{M_1} [\zeta_{N_1}] - \zeta_{N_1}(\bm y_m^1)\rp^2.\label{eq:mlvv1}
\end{align}
\end{subequations}
The variances in \eqref{eq:RV_var} are again estimated with their MC simulations \eqref{eq:mlvbound}, and $E[s_h]$ in $\widehat{\zeta}_h,\,\widehat{\zeta}_{N_1}$ is replaced by the MVR estimates $E_{M_0,M_1}[s_h]$, $\zeta_h$ and $\zeta_{N_1}$.

We would like to make two observations. First, the model and variance reduction approach described here requires $M_0$ realizations of the high-fidelity HDG output and $M_1$ realizations of the RB output, while the MC-RB approach described in the previous section requires $M$ realizations of the RB output and its error bound. If we take $N_1 = N$ then it is reasonable to consider $M_1 \approx M$. Furthermore, we take $M_0$ such that the computational cost of $M_0$ realizations of the HDG output is commensurate with that of $M$ realizations of the RB output bound. In this scenario, the two approaches have the same computational complexity. The advantage of the present approach is that it provides more accurate estimates than the MC-RB approach owing to the variance reduction. Second, unlike the MC-RB approach, the present approach does not  require {\em a posteriori} error bounds for the RB output to obtain the error bounds for our estimates of the statistical outputs. As a result, our approach  can be 
applied to problems for which {\em a posteriori} output bounds are either computationally expensive or theoretically difficult. 

\subsection{Multilevel Monte Carlo sampling}
The method can be further generalized and improved by pursuing a multilevel control variate strategy. Given $L$ different RB output models $s_{N_\ell}(\bm y), 1 \le \ell \le L,$ with $N_1 > N_2 > \ldots > N_{L}$\footnote{In our context, it is natural to number the levels from the finest RB approximation to the coarsest RB approximation because the finest RB level is closest to the HDG approximation.}, we first express the expected value as
\begin{equation*}
E[s_h] = E[s_h - s_{N_1}] + \sum_{\ell = 1}^{L-1} E [s_{N_\ell}-s_{N_{\ell+1}}]  + E[s_{N_L}]  .
\end{equation*}
We next introduce $L+1$ independent sample sets $Y^{\ell}_{M_\ell} = \{\bm y^{\ell}_m \in \Lambda,  1 \le m \le M_\ell\}$, $0 \le \ell \le L$, which are drawn in $\Lambda$ with probability density function $\rho(\bm y)$. We then define our estimate of $E[s_h]$ as
\begin{equation*}
E_{M_0,\ldots,M_L}[s_h] = E_{M_0}[s_h - s_{N_1}] +  \sum_{\ell = 1}^{L-1} E_{M_\ell} [s_{N_\ell}-s_{N_{\ell+1}}] +  E_{M_L}[s_{N_L}]  \:.
\end{equation*}
Extending the analysis in \eqref{eq:RV_exp} we apply the CLT to the multilevel case to obtain 
\begin{align*}
\lim_{M_0 \to \infty}\ldots \lim_{M_L \to \infty}& \mathrm{Pr} \left( \abs{E[s_h] - E_{M_0,\ldots,M_L}[s_h]} \le \Delta_{M_0,\ldots,M_L}^E \right) = \mathrm{erf}\lp\frac{a}{\sqrt{2}}\rp, \\
\Delta_{M_0,\ldots,M_L}^E  &=  a \sqrt{\frac{V_{M_0}[s_h - s_{N_1}]}{M_0} + \sum_{\ell=1}^{L-1}\frac{V_{M_\ell}[s_{N_\ell}-s_{N_{\ell+1}}]}{M_\ell} +\frac{V_{M_L}[s_{N_L}]}{M_L}}\:.
\end{align*}
Similarly, the estimate of the variance is defined as
\beq 
V_{M_0,\ldots,M_L}[s_h] = E_{M_0}[\zeta_h - \zeta_{N_1}] + \sum_{\ell = 1}^{L-1} E_{M_\ell} [\zeta_{N_\ell}-\zeta_{N_{\ell+1}}] +  E_{M_L}[\zeta_{N_L}]\:,
\eeq
where the auxiliary variables are $\zeta_h  := \lp s_h - E_{M_0,\ldots,M_L}[s_h] \rp^2$ and $\zeta_{N_\ell} := \lp s_{N_\ell} - E_{M_0,\ldots,M_L}[s_h] \rp^2$ for $\ell = 1,\ldots,L$. Combining the results in \eqref{eq:VarCLt}--\eqref{eq:RV_var} with the CLT leads to the following error bound for the variance estimate
\begin{align*}
\lim_{M_0 \to \infty}\ldots \lim_{M_L \to \infty}& \mathrm{Pr} \left( \abs{V[s_h] - V_{M_0,\ldots,M_L}[s_h]} \le \Delta_{M_0,\ldots,M_L}^V \right) = \mathrm{erf}\lp\frac{a}{\sqrt{2}}\rp ,\\
\Delta_{M_0,\ldots,M_L}^V  &=  a \sqrt{\frac{V_{M_0}[\zeta_h - \zeta_{N_1}]}{M_0} + \sum_{\ell=1}^{L-1}\frac{V_{M_\ell}[\zeta_{N_\ell}-\zeta_{N_{\ell+1}}]}{M_\ell}  +\frac{V_{M_L}[\zeta_{N_L}]}{M_L}} \:.
\end{align*}
Note that all expectations and variances are MC estimates through the sample sets $Y^{\ell}_{M_\ell}$ for $0 \le \ell \le L$.

We will refer to the general model and variance reduction method with a sequence of $L$ reduced basis models as the $L$-MVR method. For clarity of notation, we shall identify $s_h - s_{N_1}$ as level 0, and the subsequent $s_{N_\ell}-s_{N_{\ell+1}}$ as level $\ell$ . The method allows us to transfer the computational burden from the higher-fidelity (expensive) outputs to the lower-fidelity (inexpensive) outputs. In particular, we can choose $N_1, N_2, \ldots, N_L$ so as to have $M_0 \ll M_1 \ll \ldots \ll M_L$. Hence, the number of evaluations of the higher-fidelity outputs are significantly smaller than those of the lower-fidelity outputs, thereby resulting in a significant reduction in the overall computational cost. Finally, we address the issue of how to determine the RB dimensions $N_1, N_2, \ldots, N_L$ and the number of samples $M_0,M_1, \ldots, M_L$ to achieve a specified error tolerance and minimize the computational cost.

\subsection{Selection method}\label{sec:levsel}
Let $t_{N_\ell}$ denote the (Online) wall time to compute the RB output $s_{N_\ell}(\bm y)$ for $\ell \ge 1$,  and $t_h$ denote the wall time to compute the HDG output $s_h(\bm y)$ for any given $\bm y \in \Lambda$.  Note that $t_{N_\ell}$ depends on $N_\ell$, while $t_h$ depends on the finite element approximation spaces. The total (Online) wall time $T_L$ of the $L$-MVR and the (Online) speedup  $\pi_L$ with respect to the MC-HDG method are given by
\beql\label{eq:cost}
T_L = \lp t_h  + t_{N_1}\rp M_0 + \sum_{\ell = 1}^{L-1}  M_{\ell} \lp t_{N_\ell} +  t_{N_{\ell+1}} \rp + t_{N_L} M_L,\qquad \pi_L = \frac{t_h M}{T_L}.
\eeql
We wish to find $(N_1, N_2, \ldots, N_L)$ and $(M_0,M_1, \ldots, M_L)$  so as to minimize $T_L$, while ensuring that $\Delta_{M_0,\ldots,M_L}^E$ is equal to a specified error tolerance $\epsilon_{\rm tol}$. This error condition is satisfied if we take
\begin{equation}\label{eq:errorcondition}
a^2\frac{V_{M_L}[s_{N_L}]}{M_L} = w_L  \epsilon^2_{\rm tol},\quad  a^2\frac{V_{M_0}[s_h - s_{N_1}]}{M_0} = w_0  \epsilon^2_{\rm tol}, \quad a^2\frac{V_{M_\ell}[s_{N_\ell} - s_{N_{\ell+1}}]}{M_\ell}  = w_\ell \epsilon^2_{\rm tol}, \quad \ell \ge 1,
\end{equation}
for any given positive $w_\ell \in (0,1), \,\ell = 0,\ldots,L$ such that $w_0 + w_1 +\ldots + w_L = 1$. The choice of the weights depends on how we would like to distribute the error among the levels. We combine expressions \eqref{eq:cost}--\eqref{eq:errorcondition} to define the cost function
\beql
\label{eq:cost2}
C_L=\frac{T_L \epsilon^2_{\rm tol}}{a^2} = \frac{V_{M_0}[s_h - s_{N_1}]}{w_0} \lp t_{h}  + t_{N_1} \rp  + \sum_{\ell = 1}^{L-1} \frac{V_{M_\ell}[s_{N_{\ell}} - s_{N_{\ell+1}}]}{w_\ell} \lp t_{N_\ell} + t_{N_{\ell+1}}\rp+  t_{N_L} \frac{V_{M_L}[s_{N_L}]}{w_L} .
\eeql
We need to determine $(M_0,M_1, \ldots, M_L)$ and $(N_1, N_2, \ldots, N_L)$ that minimize $C_L$. Unfortunately, this is a nonlinear integer optimization problem which is difficult to solve exactly.  We thus solve an approximate problem as follows.

We first introduce a test sample set $Y_{\widehat{M}} = \{\widehat{\bm y}_m \in \Lambda,  1 \le m \le \widehat{M}\}$. We then precompute and store the HDG outputs $s_h(\widehat{\bm y}_m)$ for  $m = 1,\ldots,\widehat{M}$ and the RB outputs $s_{N}(\widehat{\bm y}_m)$ for  $m = 1,\ldots,\widehat{M}$ and $N = 1,\ldots,N_{\max}$.  In addition, we also precompute and store $t_h$ and $t_N$ for $N = 1, \ldots, N_{\max}$. For any given strictly decreasing $L$-tuple $\bm {I} = (I_1, I_2, \ldots, I_L) \in [1,N_{\max}]^L$ and valid weights $\bm w = (w_0,\ldots,w_L)$, we can evaluate the equivalent cost function
\beql
\label{eq:cost3}
 \widehat{C}_L(\bm {I},\bm w) = \sum_{\ell=0}^L\frac{\widehat{C}^\ell_L(\bm {I})}{w_\ell} =  \frac{V_{\widehat{M}}[s_h - s_{I_1}]}{w_0} \lp t_{h}  + t_{I_1} \rp  + \sum_{\ell = 1}^{L-1} \frac{V_{\widehat{M}}[s_{I_{\ell}} - s_{I_{\ell+1}}]}{w_\ell} \lp t_{I_\ell} + t_{I_{\ell+1}}\rp+  t_{I_L} \frac{V_{\widehat{M}}[s_{I_L}]}{w_L} ,
\eeql
with $\mcal{O}\lp(L+1)\widehat{M}\rp$ operations count, where all the variances are computed using the test sample set $Y_{\widehat{M}}$.  We now set
\begin{equation}
\label{eq:cost4}
\begin{split}
\bm N \equiv (N_1,N_2,\ldots,N_{L}) =& \arg \min_{\substack{\bm{I}}}  \:\widehat{C}_L(\bm{I},\bm{w}^{\bm{I}}),\\
&\ \textnormal{s.t.}\quad N_{\max}\geq I_1 > I_2>\ldots > I_L \geq 1
\end{split}
\end{equation}
where the weights $\bm w^{\bm{I}}$ are the minimizers of the equivalent cost for any $L$-tuple $\bm{I}$, that is
\begin{equation}
\label{eq:cost5}
\begin{split}
\bm w^{\bm I}\equiv (w_0^{\bm{I}},w_1^{\bm{I}},\ldots,w_L^{\bm{I}}) =& \arg \min_{\substack{\bm w}} \:\widehat{C}_L(\bm{I},\bm w),\\
&\ \textnormal{s.t.}\quad \sum_{\ell=0}^L w_\ell =1,\quad w_\ell > 0.
\end{split}
\end{equation}
The KKT conditions for \eqref{eq:cost5} render the optimal weights for any $L$-tuple $\bm {I}$ as
\beql\label{eq:weights}
w_\ell^{\bm{I}} = \frac{\sqrt{\widehat{C}^\ell_L(\bm {I})/\widehat{C}^0_L(\bm {I})}}{\displaystyle\sum_{\ell'=0}^L \sqrt{\widehat{C}^{\ell'}_L(\bm {I})/\widehat{C}^0_L(\bm {I})}},\qquad \ell = 0,\ldots,L.
\eeql
The minimization problem \eqref{eq:cost4} can be approximately solved, for the weights defined in \eqref{eq:weights}, by simply evaluating the cost function $\widehat{C}_L(\bm{I},\bm w^{\bm I})$ for all feasible $L$-tuples $\bm{I}$ in $\mcal{O}((L+1) \widehat{M} (N-L)\times \ldots \times(N-1)/L!)$ operations count. Even though we present here an optimal choice of the weights, any valid distribution can be employed.

Having determined the RB dimensions $\bm N$ and the weights $\bm w^{\bm N}$, we can now proceed with the MC simulations for all levels. We initially set $Y^{0}_{M_0} = Y_{\widehat{M}},$ and thus reuse $s_h(\widehat{\bm y}_m), m = 1,\ldots,\widehat{M}$. We then execute the MC processes for all the levels and enforce the error constraint  $\Delta_{M_0,\ldots,M_L}^E = \epsilon_{\rm tol}$ by adding new random parameters to the sample sets until the following inequalities 
\begin{equation*}
M_L \ge \frac{a^2V_{M_L}[s_{N_L}] }{w^{\bm N}_L\epsilon^2_{\rm tol}}, \quad M_0 \ge \frac{a^2V_{M_0}[s_h - s_{N_1}]  }{w^{\bm N}_0\epsilon_{\rm tol}},  \quad M_{\ell} \ge \frac{a^2V_{M_\ell}[s_{N_{\ell}} - s_{N_{\ell+1}}]}{ w^{\bm N}_\ell\epsilon_{\rm tol}} , \quad \ell = 1,\ldots,L-1 ,
\end{equation*}
are satisfied and the MC processes are terminated upon satisfaction of these conditions. Therefore, the sample sets $Y^{\ell}_{M_\ell}, \ell = 0,\ldots,L$ are continuously updated during the MC runs.  Finally, to provide confidence in the application of the CLT we also need to enforce that $M_\ell$ are greater than a certain threshold, say 30.

Although we have assumed that the number of levels $L$ is fixed, our approach also allows us to compare the computational costs for several values of $L$. Hence, we can determine not only the RB dimensions and the weights, but also the optimal number of levels, and it can be done efficiently evaluating expressions \eqref{eq:cost3}--\eqref{eq:cost5}. This analysis provides inexpensive means to determine the optimal multilevel structure.

\section{Numerical Results}\label{sec:results}
\subsection{A coercive example: Heat diffusion}
In the first example, we consider the one dimensional steady-state heat equation in $\mcal{D}  = (0,1)$:
\begin{subequations}\label{eq:1ddiffusion}
\begin{alignat}{2}
- \lp\kappa u_x\rp_x    & = f(x),
\quad\forall x\in\mathcal{D},\\
\kappa u_x & = 0,       \quad  \mbox{on } x=1,\\ 
u & = 0,        \quad  \mbox{on }  x=0, 
\end{alignat}
\end{subequations}
where $\kappa(x,\omega)$ is a piecewise constant function on a series of disjoint subdomains $\mcal{D}_q = ((q-1)/Q,q/Q),\;q=1,\ldots,Q$, that is, $\kappa(x,\omega) = \sum_{q=1}^Q\kappa_q(\omega)\mathbf{1}_{\mcal{D}_q}$, with $\kappa_q(\omega) \in [\gamma_q^-, \gamma_q^+]$ for all $q$. For this problem, we treat $\kappa_q(\omega)$ as i.i.d. uniform random variables; hence, we can write $\kappa(x,\omega) = \kappa(x,\sbf{y}) = \sum_{q=1}^Q y_q\mathbf{1}_{\mcal{D}_q}$, where $y_q, q = 1,\ldots,Q$ are i.i.d. random variables with uniform continuous distributions in the interval $[0.1,\, 1]$. The problem \eqref{eq:1ddiffusion} has an analytic solution given by
\beql\label{eq:exactsolution}
u(x,\sbf{y}) = \int_0^x \lp \frac{1}{\kappa(z,\sbf{y})}\int_z^1 f(\xi)d\xi \rp dz.
\eeql
The observable quantity is the average temperature on the domain, namely, $s(\sbf{y}) = \int_0^1 u(x,\sbf{y})$. The output $s$, its expectation $E[s]$ and variance $V[s]$ have closed analytic forms, thereby $s$ will be here used instead of the HDG output $s_h$. Numerical results for a constant source term $f(x) = 1$ and $Q = 10$ are presented below.

 \begin{figure}[h!]
\centering
\subcaptionbox{Exact solution $u(x,\sbf{y})$ for several realizations of the diffusivity field.\label{fig:diffusivity}}
[7.5cm]{\includegraphics[scale = 0.45]{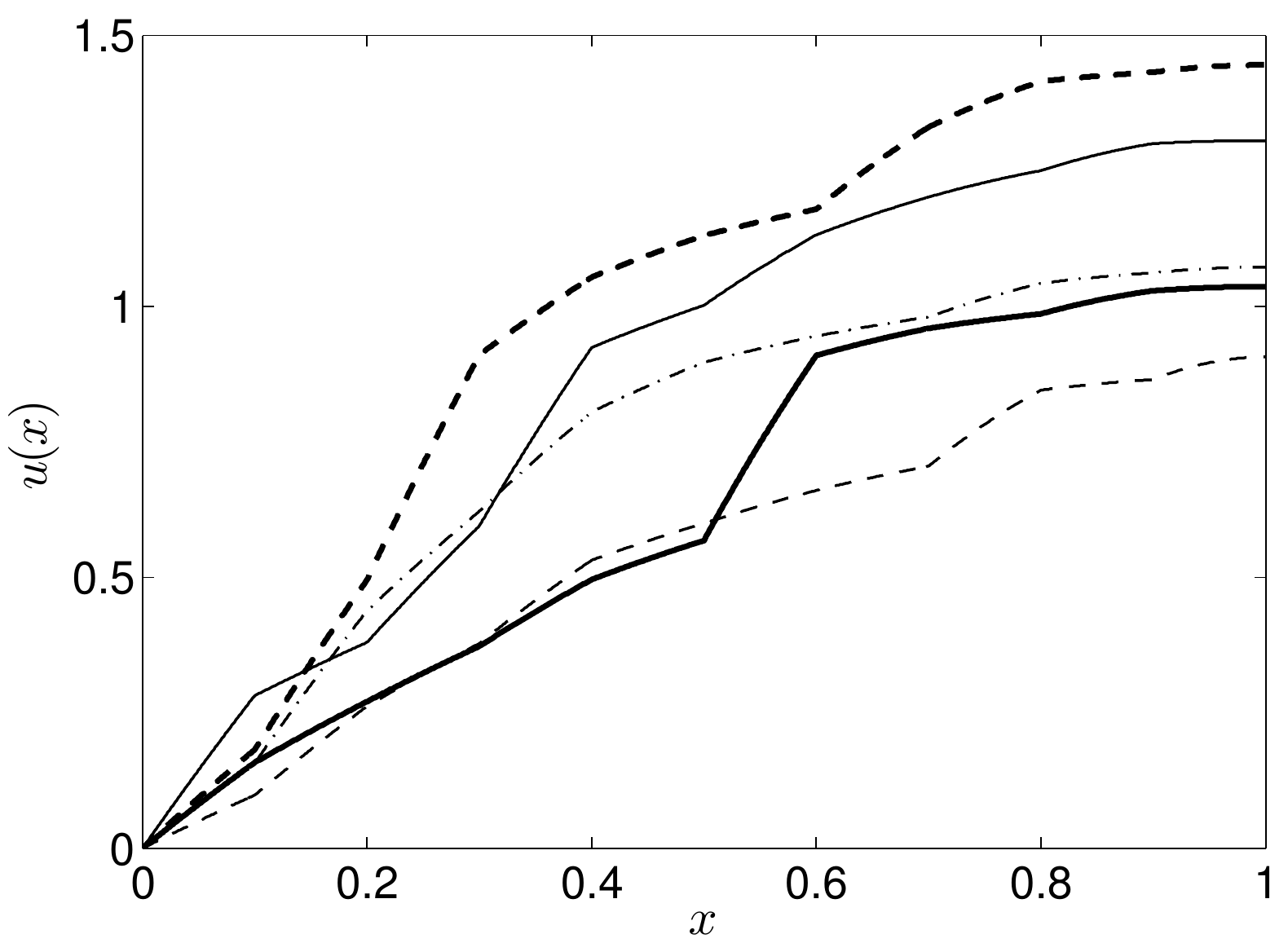}}
\hfill\subcaptionbox{Average output error bound $\Delta^s_{N,\rm avg}$ and average output error $\epsilon_{N,\rm avg}$ for $\sbf{y}\in Y_{\widehat{M}}$ vs. RB size N. \label{fig:rb1d}}
[7.5cm]{\includegraphics[scale = 0.45]{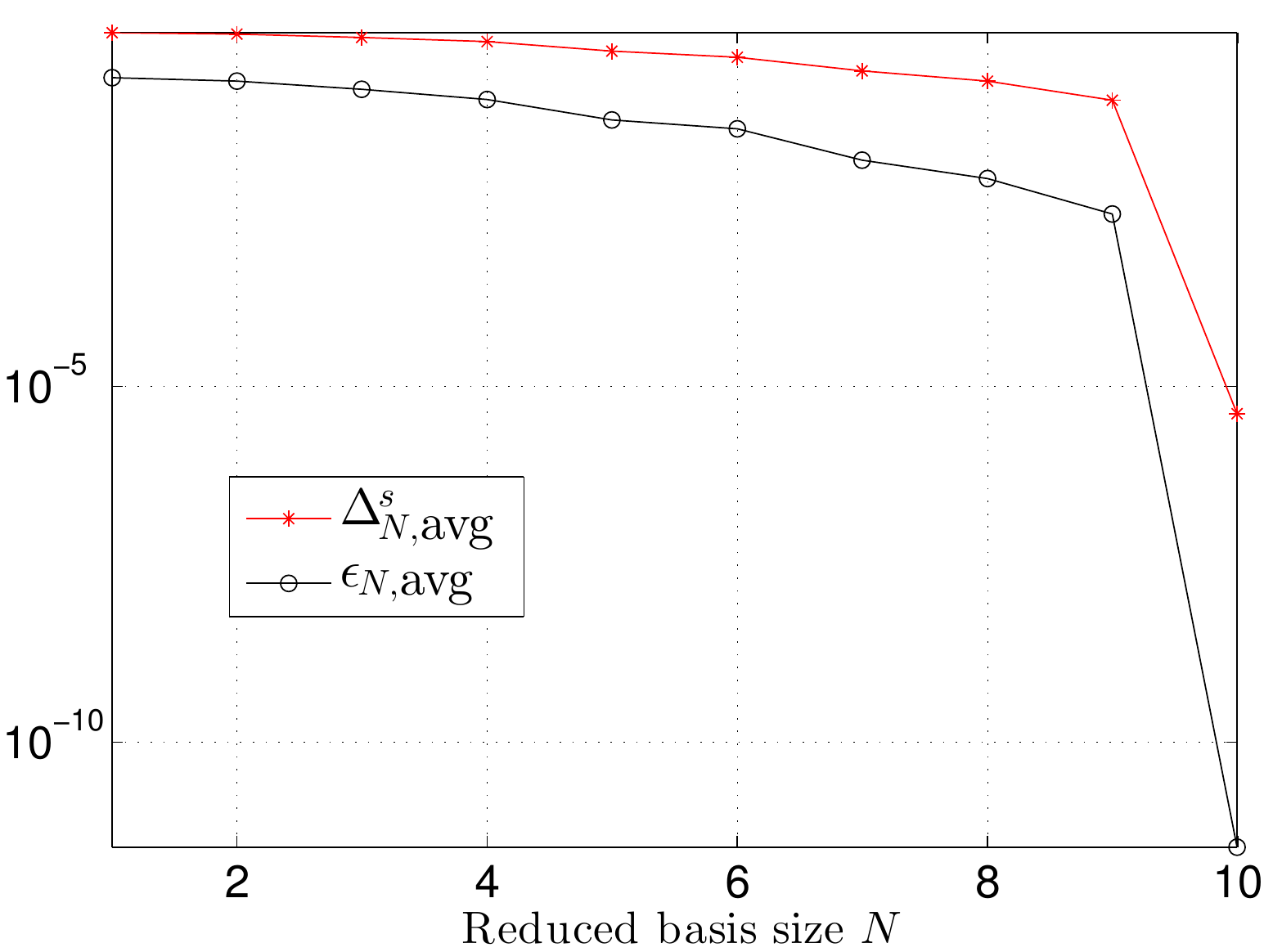}}
\caption{Representative solutions and the RB convergence of the heat diffusion example.}
\end{figure}

\begin{table}[htdp]
\footnotesize
\begin{center}
\begin{tabular}{|c||c|c|c|c|c|c|c|c|}
\hline
  & \multicolumn{2}{|c}{MC-HDG} & \multicolumn{2}{|c}{MC-RB ($N=9$)} &  \multicolumn{2}{|c|}{$1$-MVR ($N_1 = 5$)} \\
\hline
$M$  & $\abs{E[s] - E_M[s]}$ & $\Delta^E_{h,M}$ & $\abs{E[s] - E_M[s_N]}$ &$\widetilde{\Delta}^E_{N,M}$ & $\abs{E[s] - E_{M_0,M_1}[s]}$ & $\Delta^E_{M_0,M_1}$  \\
\hline
\hline
$10^2$  &  $2.10\textnormal{e}{-}2$ & $2.6\textnormal{e}{-}2\:a$ &  $2.19\textnormal{e}{-}2$ &  $1.1\textnormal{e}{-}1 + 7.0\textnormal{e}{-}2\:a$  & $2.30\textnormal{e}{-}2$&   $2.9\textnormal{e}{-}2\:a $ \\
\hline
{$10^3$} &  $6.59\textnormal{e}{-}2$ & $8.3\textnormal{e}{-}3\:a$&  $7.06\textnormal{e}{-}3$ & $1.1\textnormal{e}{-}1 + 2.2\textnormal{e}{-}2\:a$  & $7.30\textnormal{e}{-}3$ &   $9.2\textnormal{e}{-}3\:a $ \\
\hline
{$10^4$} &  $2.08\textnormal{e}{-}3$ & $2.6\textnormal{e}{-}3\:a$&  $3.21\textnormal{e}{-}3$ & $1.1\textnormal{e}{-}1 + 7.0\textnormal{e}{-}3\:a$ & $2.32\textnormal{e}{-}3$&  $2.9\textnormal{e}{-}3\:a $ \\
\hline
{$10^5$} &  $6.57\textnormal{e}{-}4$ & $8.3\textnormal{e}{-}4\:a$&  $2.88\textnormal{e}{-}3$ & $1.1\textnormal{e}{-}1 + 2.2\textnormal{e}{-}3\:a$  & $7.29\textnormal{e}{-}4$ &   $9.2\textnormal{e}{-}4\:a $ \\
\hline
{$10^6$} &  $2.11\textnormal{e}{-}4$ & $2.6\textnormal{e}{-}4\:a$&  $2.90{e}{-}3$ & $1.1\textnormal{e}{-}1 + 7.0\textnormal{e}{-}4\:a$ & $2.32\textnormal{e}{-}4$&  $2.9\textnormal{e}{-}4\:a $ \\
\hline
\end{tabular}
\end{center}
\caption{The expectation error and its error bound for different values of $M$ for the MC-HDG method, the MC-RB method and the $1$-MVR method.}
\label{tab:1derror}
\end{table}%
We show in Figure \ref{fig:diffusivity} different realizations of the exact solution in \eqref{eq:exactsolution}. Since the output is compliant, the dual problem coincides with the primal problem. We thus need to construct the reduced basis approximation for the primal problem only. Furthermore, since the bilinear form is coercive and the parameters are positive, we compute the stability constant ${\beta}_h(\sbf{y})$ using a bound conditioner technique \cite{veroy2002posteriori}, which greatly simplifies the process. Our reduced basis is constructed with $N_{\max} = 10$.  We show in Figure 1(b) the average output error $\epsilon_{N,\rm avg}$ and the average output error bound $\Delta^s_{N, \rm avg}$ as a function of $N$, where $\epsilon_{N,\rm avg} = \sum_{\bm y \in Y_{\widehat{M}}} \abs{s(\sbf{y})- s_N(\sbf{y})}/\widehat{M} $ and $\Delta^s_{N, \rm avg} = \sum_{\bm y \in Y_{\widehat{M}}}  \Delta^s_{N}(\bm y)/\widehat{M} $, being $Y_{\widehat{M}}$ a test set of $\widehat{M} = 1000$ samples. We observe 
that the average output error and the 
average output error bound converge slowly up to $N = 9$ and drop rapidly at $N = 10$. This is because of the nature of the particular problem which requires $N = Q$ basis functions to capture all the possible solutions. When we use $N < Q$, we do not have enough basis functions to represent all the possible solutions, which in turn causes a slow convergence of the reduced basis approximation.
\begin{table}[htdp]
\footnotesize
\begin{center}
\begin{tabular}{|c||c|c|c|c|c|c|}
\hline
  & \multicolumn{2}{|c}{MC-HDG} & \multicolumn{2}{|c|}{$1$-MVR ($N_1 = 5$)} \\
\hline
$M$  & $\abs{V[s] - V_M[s]}$ & $\Delta^V_{h,M}$ & $\abs{V[s] - V_{M_0,M_1}[s]}$ & $\Delta^V_{M_0,M_1}$  \\
\hline
\hline
$10^2$  &  $9.21\textnormal{e}{-}3$ & $1.1\textnormal{e}{-}2\:a$ & $1.24\textnormal{e}{-}2$&   $1.5\textnormal{e}{-}2\:a $ \\
\hline
{$10^3$} &  $2.91\textnormal{e}{-}3$ & $3.6\textnormal{e}{-}3\:a$& $4.06\textnormal{e}{-}3$ &   $5.0\textnormal{e}{-}3\:a $ \\
\hline
{$10^4$} &  $9.20\textnormal{e}{-}4$ & $1.2\textnormal{e}{-}3\:a$& $1.29\textnormal{e}{-}3$&  $1.6\textnormal{e}{-}3\:a $ \\
\hline
{$10^5$} &  $2.91\textnormal{e}{-}4$ & $3.6\textnormal{e}{-}4\:a$& $4.07\textnormal{e}{-}4$ &   $5.1\textnormal{e}{-}4\:a $ \\
\hline
{$10^6$} &  $9.29\textnormal{e}{-}5$ & $1.1\textnormal{e}{-}4\:a$& $1.27\textnormal{e}{-}4$&  $1.6\textnormal{e}{-}4\:a $ \\
\hline
\end{tabular}
\end{center}
\caption{The variance error and its error bound for different values of $M$ for the MC-HDG method and the $1$-MVR method.}
\label{tab:1derror2}
\end{table}%

We now compare the performance of the MC-HDG  method with a uniform mesh of $h = 1/10$, the MC-RB with $N = 9$, and the $L$-MVR with $L = 1$, $N_1 = 5$, $M_1 = M$ and $M_0 = M/10$  in estimating $E[s]$ and $V[s]$ as a function of $M$. For each $M$ value we repeat the simulations $H=1000$ times, and present in Tables \ref{tab:1derror} and \ref{tab:1derror2} the average values of the absolute errors and error bounds for the expectation and the variance respectively. We observe that the $1$-MVR method significantly outperforms the MC-RB method. The improvement is noticeable when we increase $M$, since the MC-RB method stagnates around $2.9\times 10^{-3}$ whereas $1$-MVR keeps reducing the error as the square root of the number of samples. The stagnation is caused by the inherent bias arising from the reduced basis method with $N=9$, which provides outputs with a level of error of $3 \times 10^{-3}$, as seen in Figure 1(b). With the MC-RB method we are unable to achieve more accurate estimators than the precision of the reduced basis output. Furthermore, 
the error bound of the MC-RB method,  $\widetilde{\Delta}_{N,M}^{E}$ as defined in (\ref{eq:CLTMCRBbound}), is the sum of two terms: the first term ${\Delta}_{N,M}^{E}$ does not depend on the confidence level, whereas the second  term $a\sqrt{{(V_M[s_N] + \Delta_{N,M}^{V} )}/{M}}$ does. We observe from Table \ref{tab:1derror} that increasing $M$ does not improve the error bound of the MC-RB method since it is dominated by ${\Delta}_{N,M}^{E}$, which in many cases can be overly pessimistic. Variance estimations for the MC-RB are not included in Table \ref{tab:1derror2}, as they can only be worse than the expectation results.

The $1$-MVR method does not suffer from this stagnation owing to the fact that it directly approximates $E[s]$ instead of $E[s_N]$. As a result, the expectation and variance error of the $1$-MVR method can be made arbitrarily small. The same behavior is observed in the error bounds $\Delta^E_{M_0,M_1}$ and $\Delta^V_{M_0,M_1}$ defined in \eqref{eq:mlebound} and \eqref{eq:mlvbound}, which agree with the Monte Carlo dependence on the square root of the sample size. Even though the accuracy of the estimators and the 
sharpness of the bounds for 1-MVR is slightly worse than that of MC-HDG, the former performs ten times less full model evaluations than the latter. These numerical results show a considerable gain for model and variance reduction.

\subsection{A noncoercive example: Acoustic wave propagation}

\begin{figure}[ht!]
\centering
\subcaptionbox{Geometry of the wave propagation problem. The source generates a wave that propagates through the medium. \label{fig:wavegeometry}}
[7.4cm]{\includegraphics[scale = .9]{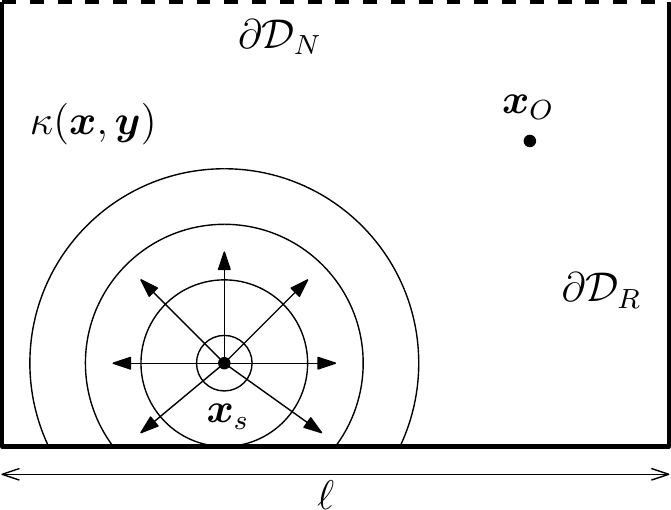}}
\hfill\subcaptionbox{Triangular mesh of 1420 elements for the wave problem. Higher resolution is appreciated at the source and output locations. \label{fig:meshwave}}
[7.4cm]{\includegraphics[scale = 0.6]{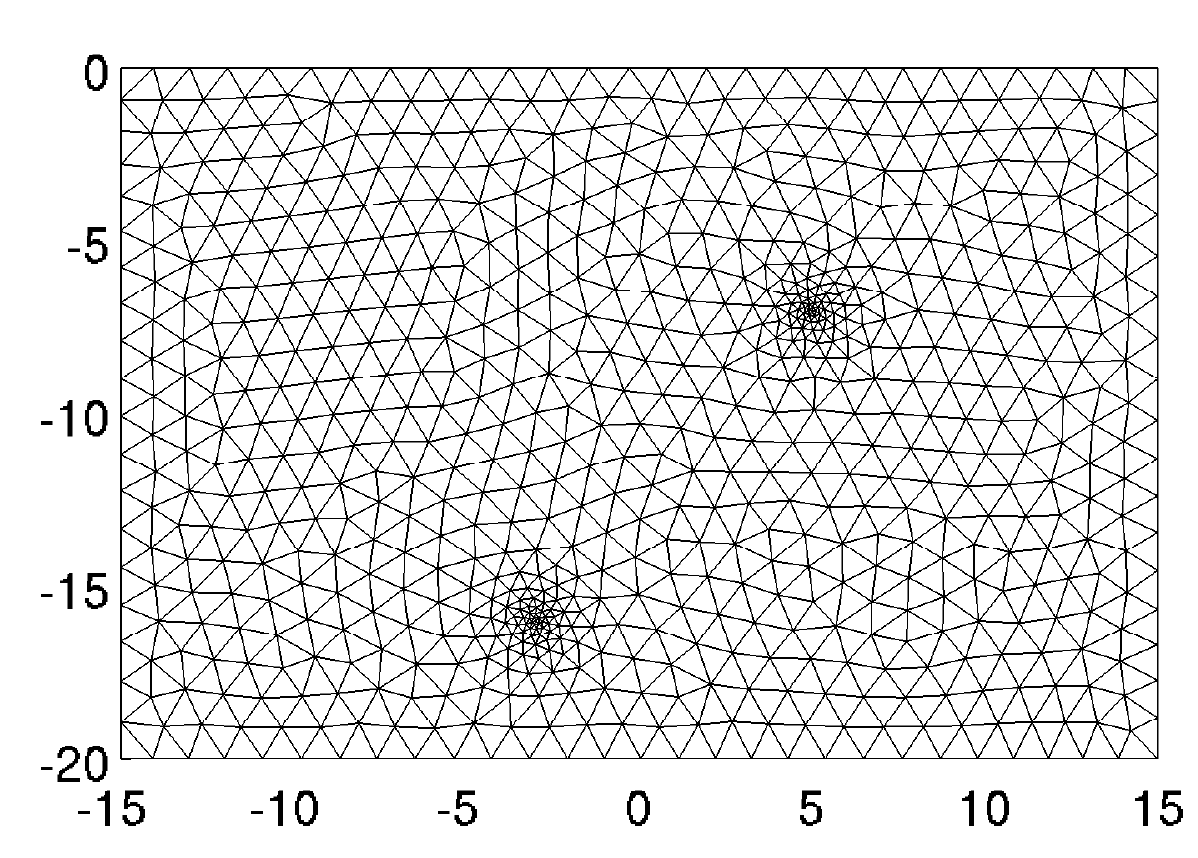}}
\hfill\subcaptionbox{Random realization of the wave amplitude field, using elements of order $p=4$ for a total of $\mcal{N} = 32180$ degrees of freedom. \label{fig:wavesolution}}
[7.4cm]{\includegraphics[scale = 0.5]{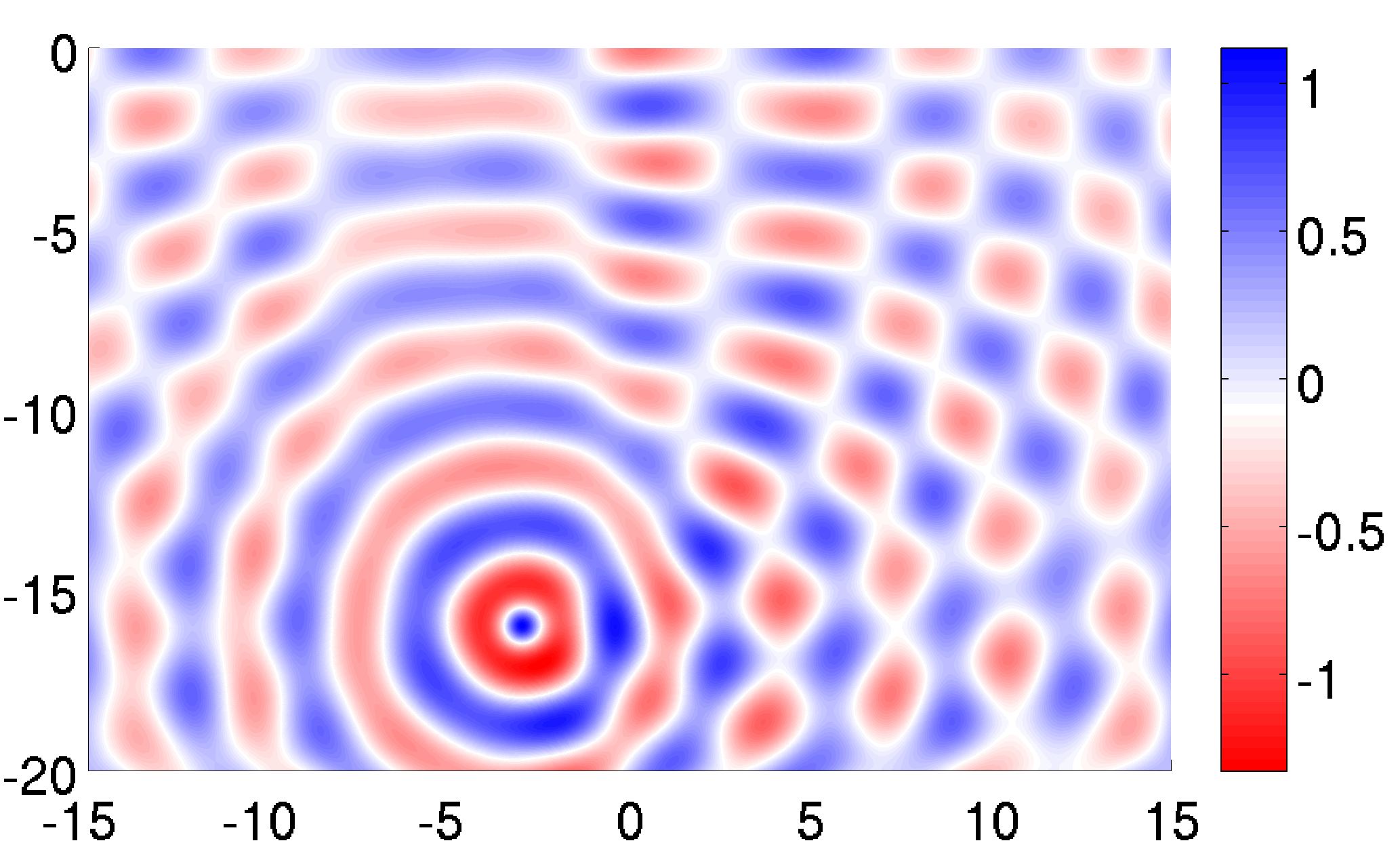}}
\hfill\subcaptionbox{Average output error bound $\Delta^s_{N,\rm avg}$ and average output error $\epsilon_{N,\rm avg}$ for $\sbf{y}\in Y_{\widehat{M}}$ vs. RB size N.\label{fig:offlineConv}}
[7.4cm]{\includegraphics[scale = 0.5]{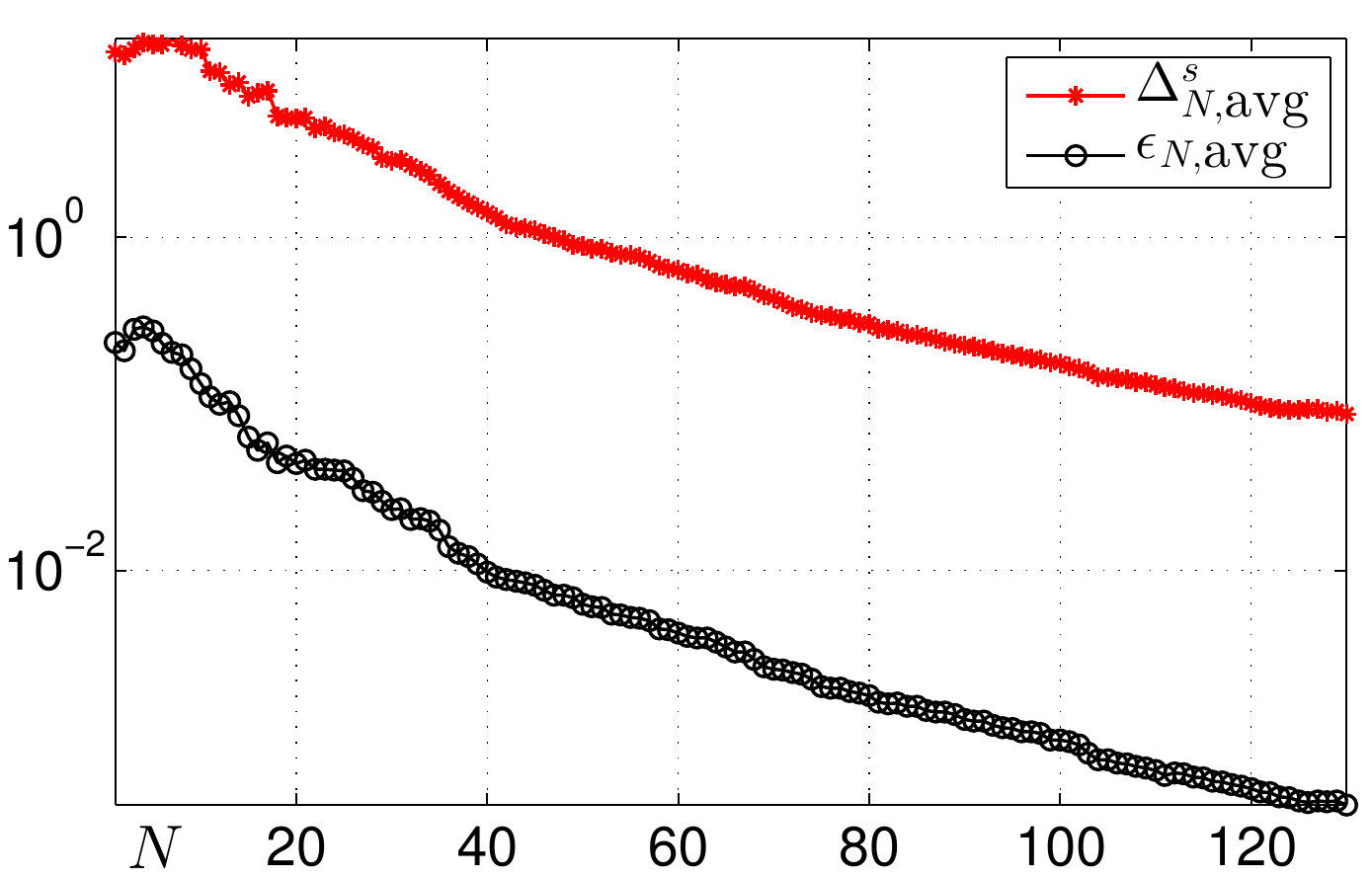}}
\caption{Problem specification, representative solution, and RB convergence of the wave propagation example}
\end{figure}

We consider a wave propagation problem as depicted in Figure \ref{fig:wavegeometry}. A wave is excited by a Gaussian source term $f $ centered at $\bx_s$ and propagates through a heterogeneous medium $\kappa(\bx,\bm y)$. The governing equation for this model problem is given by
\begin{alignat*}{2}
-\nabla\cdot (\kappa\ggrad u)  - k^2 u&= f,\qquad& \forall \bx& \in  \mathcal{D}\;, \\
 \kappa \ggrad u \cdot \bn - i k u\ &= 0,\qquad& \forall \bx& \in \partial \mathcal{D}_R \;,  \\
  \kappa \ggrad u \cdot \bn  &= 0,\qquad& \forall \bx& \in \partial \mathcal{D}_N \;,
\end{alignat*}
where $f = \dfrac{10}{\sqrt{2\pi}\sigma_s}\exp \left(- \dfrac{(x_1-x_{s1})^2+(x_2-x_{s2})^2}{2\sigma_s^2} \right)$ for $\boldsymbol{x}_s = (-3,-16)$ and $\sigma_s = 0.25$ is the source term and $k = \sqrt{2}$ is the wavenumber. Here the physical domain is $\mcal{D} = [-15,15]\times[-20,0]$. To describe the $\kappa(\bx,\sbf{y})$ field we use the example described in \cite{chen2014comparison}, namely
\beq
\kappa(\bx,\sbf{y}) = \overline{\kappa} + \sigma y_1\sqrt{\frac{\lambda_0}{2}} + \sigma\sum_{n=1}^{8} \sqrt{\lambda_n}\lp \sin \lp n\pi \frac{x_1+\ell/2}{\ell}\rp y_{2n} + \cos \lp n\pi \frac{x_1+\ell/2}{\ell}\rp y_{2n+1}  \rp,  
\eeq
where
\beq
\sqrt{\lambda_n} = \lp \sqrt{\pi} L_c\rp^{1/2}\exp \lp -\frac{(n\pi L_c)^2}{8} \rp,\quad n = 0,\ldots,8,
\eeq
for $\overline{\kappa} = 1,\sigma = 1/10$, and $L_c= 1/12$. Here the random variables $y_n$ for $n = 1,\ldots,Q=17$ are uncorrelated and uniformly distributed with zero mean and unit variance. Hence, we can write $\kappa(\bm x, \bm y)$ in the form of the affine expansion (\ref{affinekappa}) with $\Lambda = [-\sqrt{3},\sqrt{3}]^Q$. We consider the following output  
\begin{equation*}
s(\bm y) = \dfrac{1}{\sqrt{2\pi}\sigma_O}\int_{\mathcal D} \Re(u(\bm y)) \exp \left(- \frac{(x_1-x_{O1})^2+(x_2-x_{O2})^2}{2\sigma_O^2} \right) d \bm x, 
\end{equation*}
for $\bm x_O =  (5,-7)$ and $\sigma_O = 0.25$, which corresponds to the real part of the amplitude at $\bm x_O$ regularized by a Gaussian field. The physical domain is discretized into a triangular mesh of 1420 elements as shown in Figure \ref{fig:meshwave} and polynomials of degree $p = 4$ are used to represent the numerical solution $u_h(\bm y)$. Figure \ref{fig:wavesolution} depicts a realization of the numerical solution obtained using the HDG method. 

Since the exact values of $E[s_h]$ and $V[s_h]$ are not known, we approximate them using the MC-HDG method with a random sample set $Y_{M^*}$ of $M^* = 6.5\times 10^7$ and obtain $E_{M^*}[s_h] = 0.2576$ and $V_{M^*}[s_h] = 0.0596$ for a statistical error of  $10^{-4}$ and $3.1 \times 10^{-5}$ respectively, corresponding to $0.999$ confidence level ($a = 3.3$). We are going to use $E[s_h]=E_{M^*}[s_h]$ and $V[s_h]=V_{M^*}[s_h]$ as the reference values to evaluate the performance of our method. We consider an error tolerance of $10^{-3}$ and a confidence level $0.95$ for our estimators. To achieve this level of accuracy, the MC-HDG method requires a random sample set $Y_M$ of size $M = 238447$ to compute the MC-HDG estimators $E_{M}[s_h],\,V_{M}[s_h]$.

\begin{figure}[h!]
\centering
\subcaptionbox{Speedup $\pi_1$ and cost functions $C_1,\,\widehat{C}_1$.\label{fig:1-speedup}}
[4.92cm]{\includegraphics[scale = 0.31]{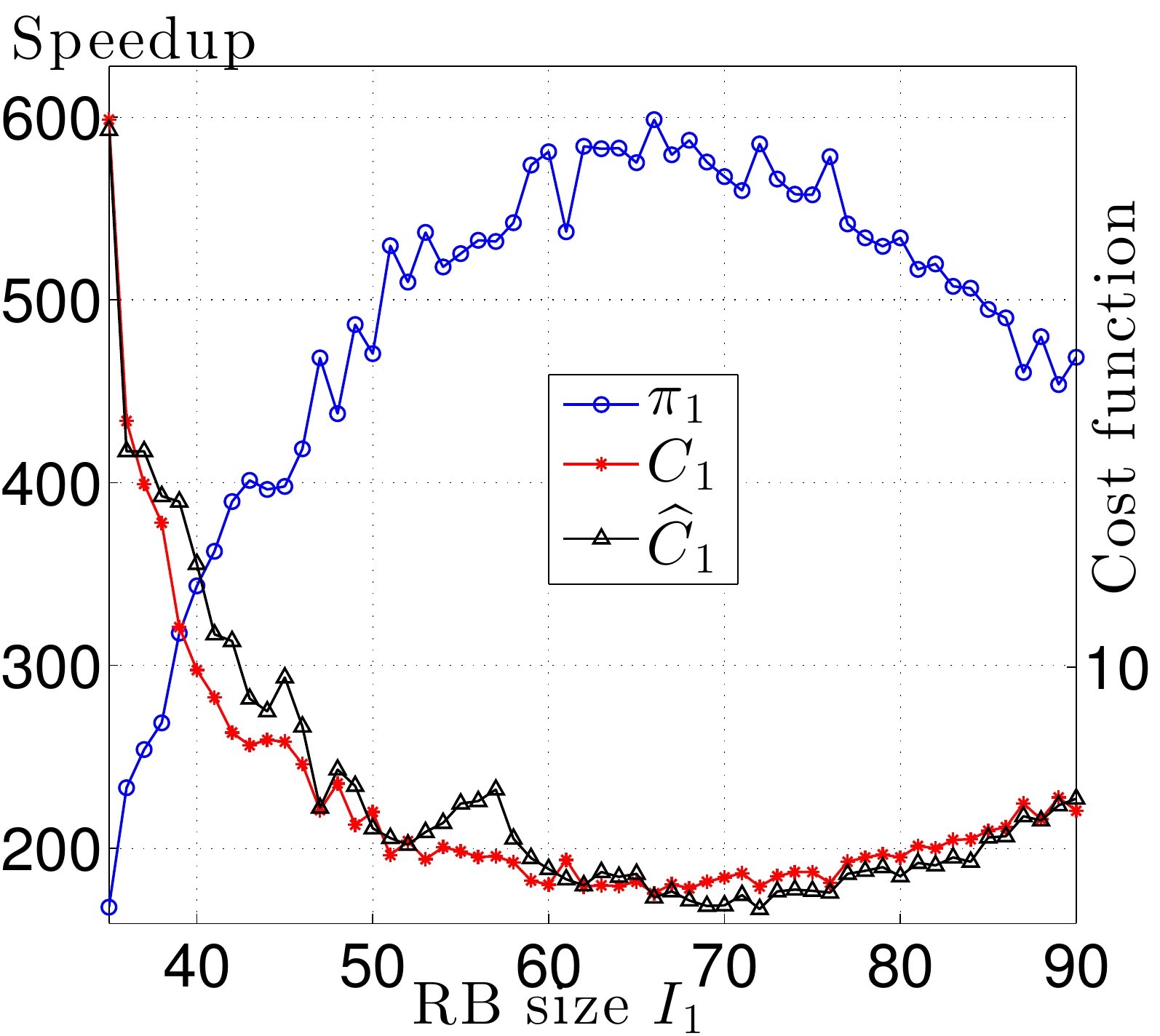}}
\hfill\subcaptionbox{Variance of $s_h - s_{I_1}$ and relative time to solve RB output. \label{fig:1-mvr}}
[4.92cm]{\includegraphics[scale = 0.31]{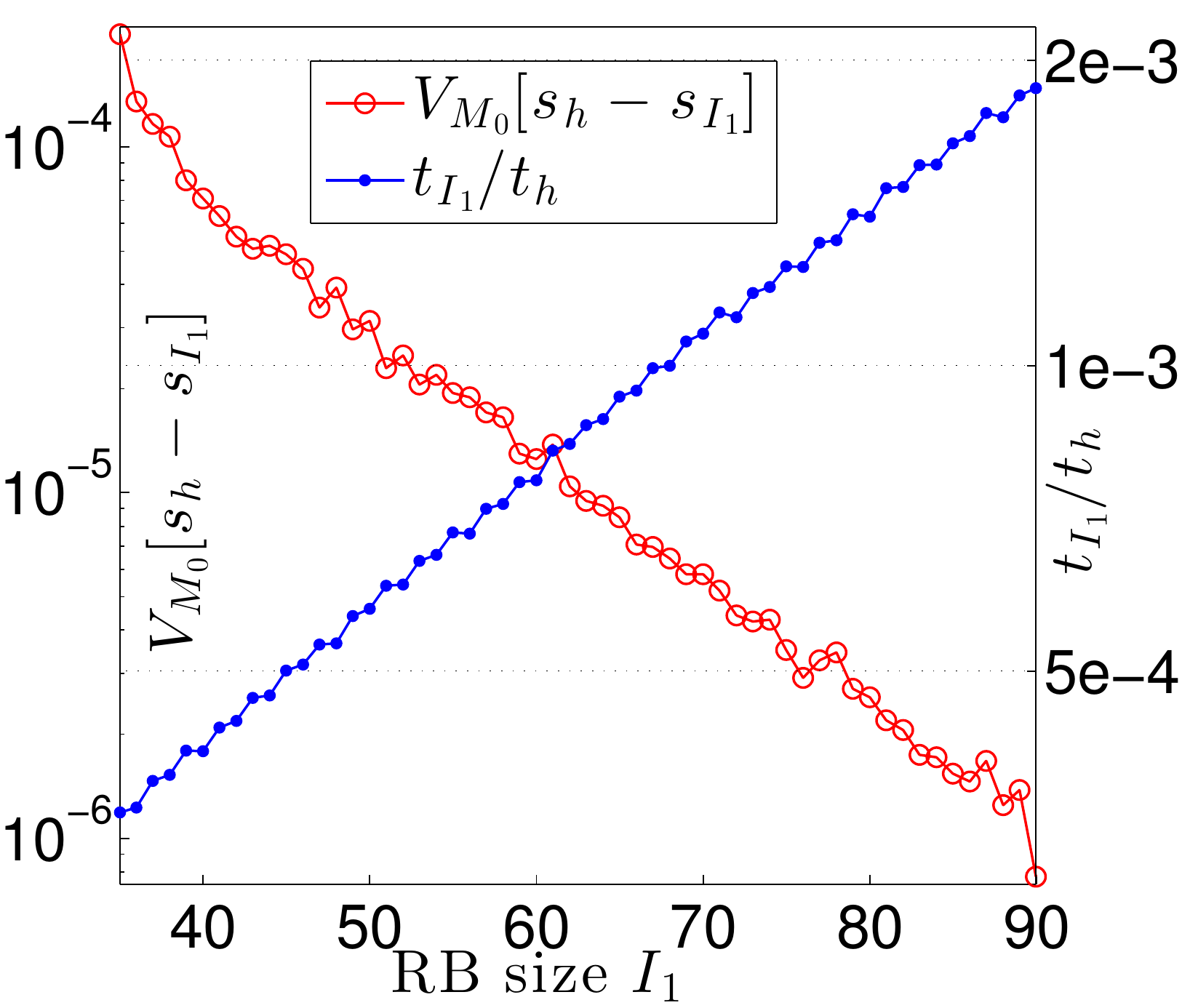}}
\hfill\subcaptionbox{Optimal weights.\label{fig:1-weights}}
[4.92cm]{\includegraphics[scale = 0.31]{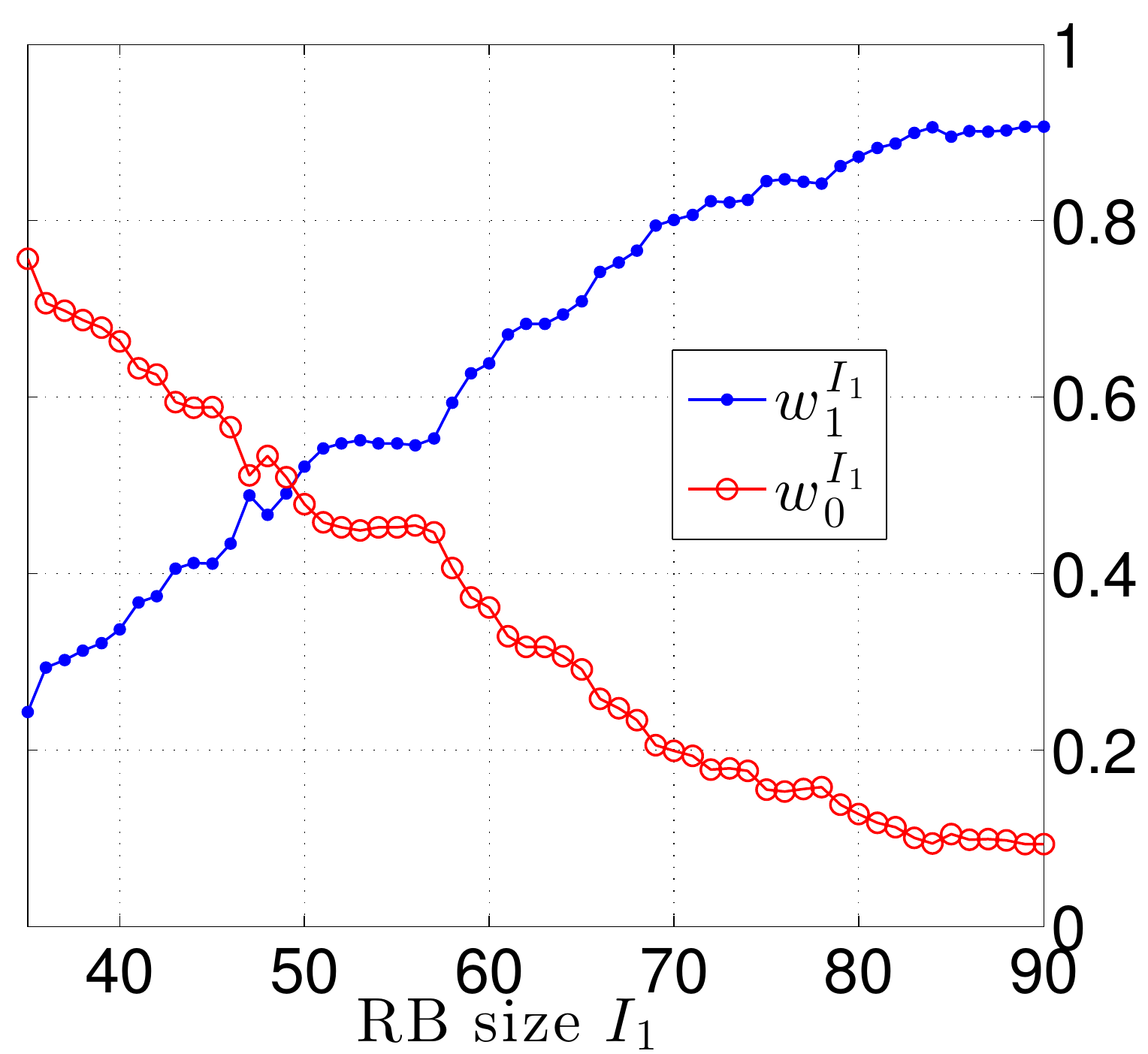}}
   \hfill\subcaptionbox{Expectation (left) and variance (right) estimators with 95\% confidence interval.\label{fig:errorbar}}
  [15cm]{\includegraphics[scale = 0.38]{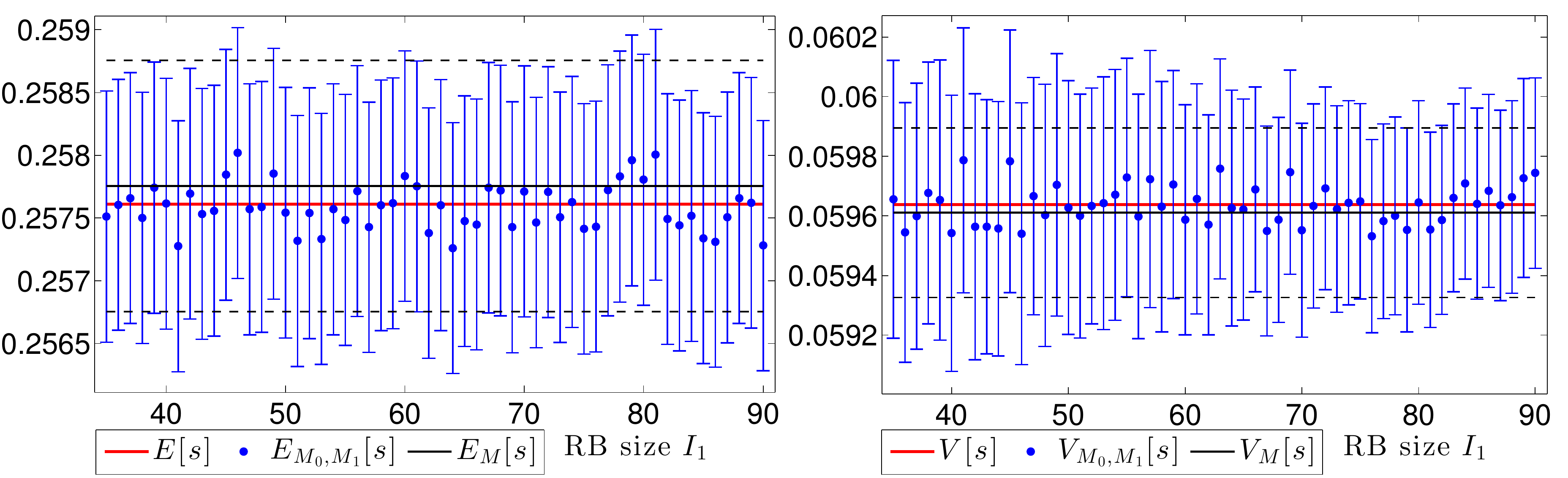}}
\caption{Result for the $1$-MVR method vs RB size $I_1$.}
\end{figure}

We next pursue the RB method and show in Figure \ref{fig:offlineConv} the average output error $\epsilon_{N,\rm avg}$ and its error bound $\Delta^s_{N,\rm avg}$ as a function of $N$. Here $\epsilon_{N,\rm avg} = \sum_{\bm y \in Y_{\widehat{M}}} \abs{s(\sbf{y})- s_N(\sbf{y})}/\widehat{M}$ and $\Delta^s_{N, \rm avg} = \sum_{\bm y \in Y_{\widehat{M}}}  \Delta^s_{N}(\bm y)/\widehat{M}$, where $Y_{\widehat{M}}$ is a test set of $\widehat{M} = 100$ samples. We observe that the RB error bound is about two orders of magnitude larger than the output error. The slow convergence rate of the RB error bound is expected because the problem is non coercive and has (many) $Q=17$ parameters.  Since the RB error bounds are quite pessimistic, we will not consider the MC-RB method to compute the statistical outputs and their error bounds.
\begin{figure}[h!]
 \centering
 \subcaptionbox{Computational speedup $\pi_2$.\label{fig:2-speedup}}
 [4.92cm]{\includegraphics[scale = .28]{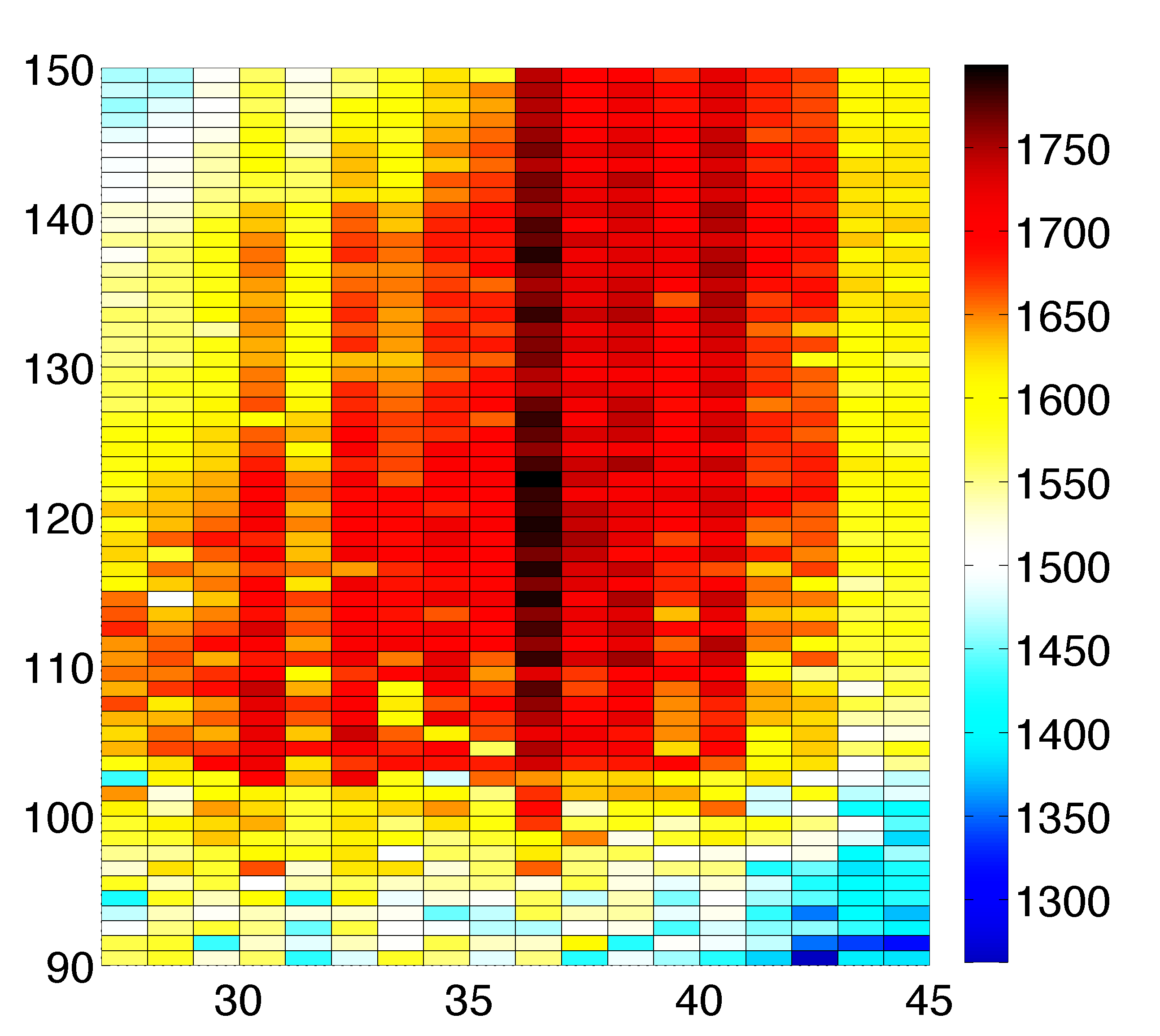}}
 \hfill \subcaptionbox{Cost function $C_2$.\label{fig:2-cost}}
  [4.92cm]{\includegraphics[scale = .28]{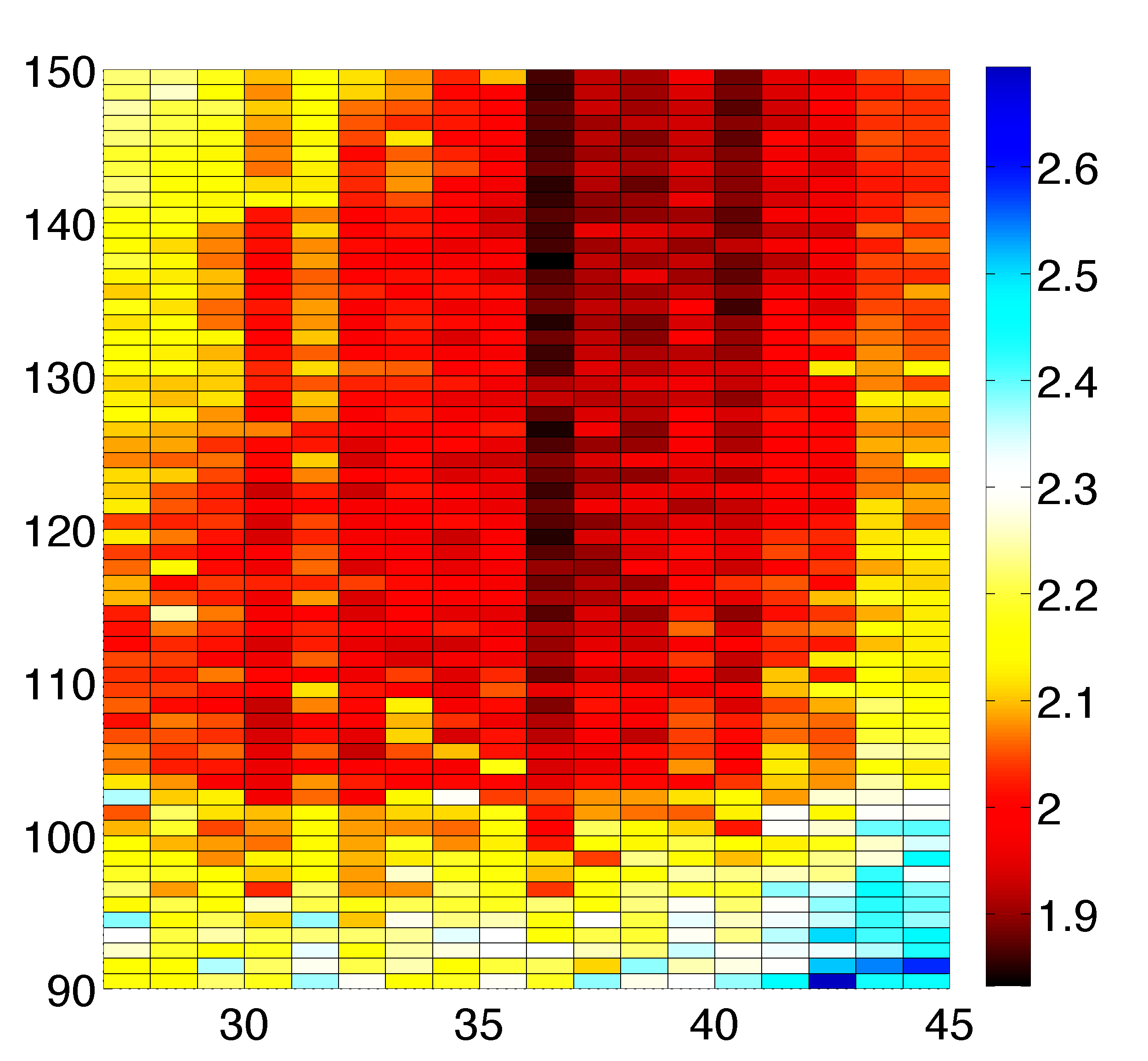}}
  \hfill \subcaptionbox{Equivalent cost function $\widehat{C}_2$.\label{fig:2-surrcost}}
   [4.92cm]{\includegraphics[scale = .28]{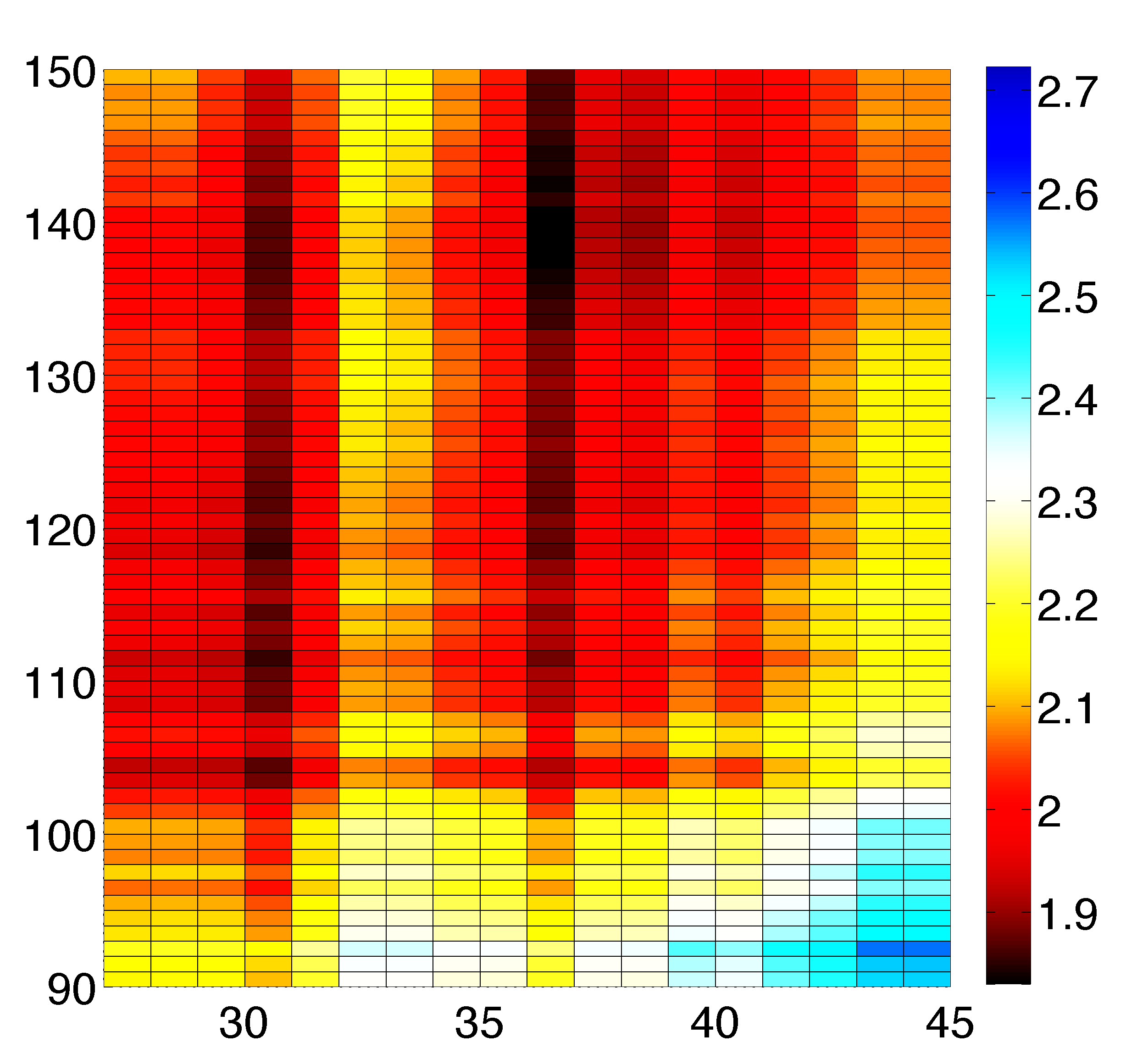}}
   \hfill 
  \subcaptionbox{Expectation (left) and variance (right) estimators with 95\% confidence interval.\label{fig:2-accuracy}}
  [15cm]{\includegraphics[scale = 0.3]{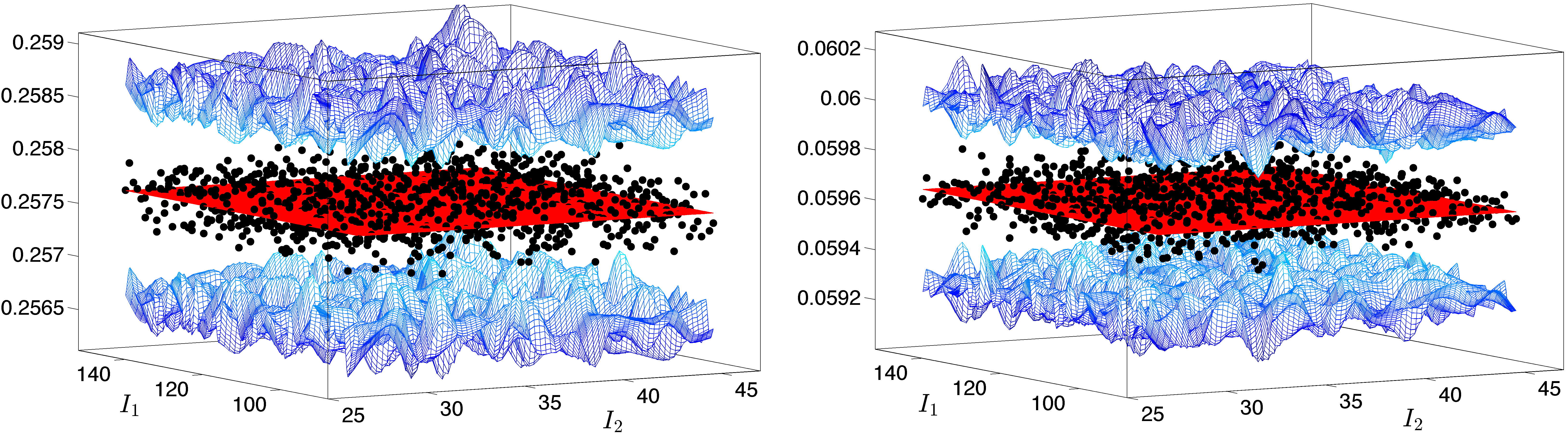}}
  \caption{Results for the $2$-MVR method vs RB sizes $I_2$ and $I_1$.}
 \end{figure}

We now turn to the $1$-MVR and enforce the tolerance $\Delta_{M_0,M_1}^E = \epsilon_{\rm tol} = 10^{-3}$ and the confidence level of 0.95. We depict in Figure \ref{fig:1-speedup} the computational speedup $\pi_1$ relative to the MC-HDG method, the original cost function $C_1$, and the equivalent cost function $\widehat{C}_1$ as a function of the RB dimension $I_1$-- results for each level size are averaged 8 times. The equivalent cost function $\widehat{C}_1$ approximates the original cost function $C_1$ reasonably well, despite being drastically less expensive to evaluate than the true cost function (and available \emph{a priori}). The equivalent cost is minimized at $N_1 = 72$, requiring $(M_0,M_1) = (96,278684)$ for a speedup $\pi_1 = 585$, whereas the true cost yields an optimal RB dimension $N_1 = 66$ and sample sizes $(M_0,M_1) = (106,308836)$ that achieve a speedup $\pi_1 = 599$. The model and variance reduction strategy is represented in Figures \ref{fig:1-mvr}-\ref{fig:1-weights}. For small reduced 
basis size more accuracy is demanded for level 1, that is $w_1 < w_0$, since the RB model is very inexpensive to evaluate and the variance of  $s_h-s_{I_1}$  is large. Conversely, for increasing $I_1$ the RB model becomes more costly to compute, whereas the variance of $s_h-s_{I_1}$ rapidly decreases -- therefore requiring very few full model evaluations. This change of behavior is detected by the level selection method by setting $w_0 <w_1$, that is requiring higher accuracy for level 0. The model and variance reduction method therefore seeks a balance between these two phenomena to achieve optimal efficiency. We next show in Figure \ref{fig:errorbar} the $1$-MVR expectation and variance as well as their error bars as a function of $I_1$. The expectation error bound is equal to the prescribed tolerance $\epsilon_{\rm tol} = 10^{-3}$, while the error bound for the variance $\Delta_{M_0,M_1}^{V}$ decreases from about $4.5 \times 10^{-4}$ for small $I_1$ to about $3.2 \times 10^{-4}$ for larger reduced basis 
size.

We proceed analogously for the $2$-MVR method computing the speedup with respect to the MC-HDG method for several values of level sizes $I_1$ and $I_2$. The computational gain is presented in Figure \ref{fig:2-speedup}, and the optimum speedup $\pi_2 = 1800$ is reached for $(N_1,N_2) = (122,36)$, which needs $(M_0,M_1,M_2) = (30,4102,286152)$.  The true cost and the equivalent cost present a very similar behavior, compare Figures \ref{fig:2-cost}-\ref{fig:2-surrcost}, and the level selection method recovers as low fidelity models the bases $(N_1,N_2) = (139,36)$, for sample sizes $(M_0,M_1,M_2) = (30,3483,281824)$ and a speedup $\pi_2 = 1778$. The accuracy of the estimators is shown in Figure \ref{fig:2-accuracy}, and we observe that the true moments lie in all cases inside the confidence intervals (displayed as a surface) computed by the model and variance reduction method.
 \begin{table}[ht]
 \footnotesize
 \centering
 \begin{tabular}{|l||l|l|l|l|l|}
 \hline
 &  $L=1$ & $L=2$ & $L=3$ & $L=4$ &$L=5$\\
 \hline\hline
 $( N_1,\ldots,N_L)$ & 72 & $(139,36)$ & $(150,52,16)$ & $(150,84,47,16)$& $(150,118,84,47,16)$\\ \hline
 $\pi_L$ & 585 & 1778 & 2277 & 2363 & 2354\\ \hline
 $C_L/C_4$ & 4.03 & 1.33 & 1.04 & 1 & 1.00\\ \hline
 $\widehat{C}_L/\widehat{C}_4$  & 4.37 & 1.38 & 1.08 &1  &1.00\\ \hline
 \end{tabular}
 \caption{Predicted optimal level sizes, speedup $\pi_L$ with respect to MC-HDG cost and relative real and equivalent multilevel cost.}\label{tab:multilevel}
 \end{table}
 \begin{figure}[h!]
  \centering
  \subcaptionbox{Distribution of weights.\label{fig:L-weights}}
     [4.92cm]{\includegraphics[width = 4.85cm]{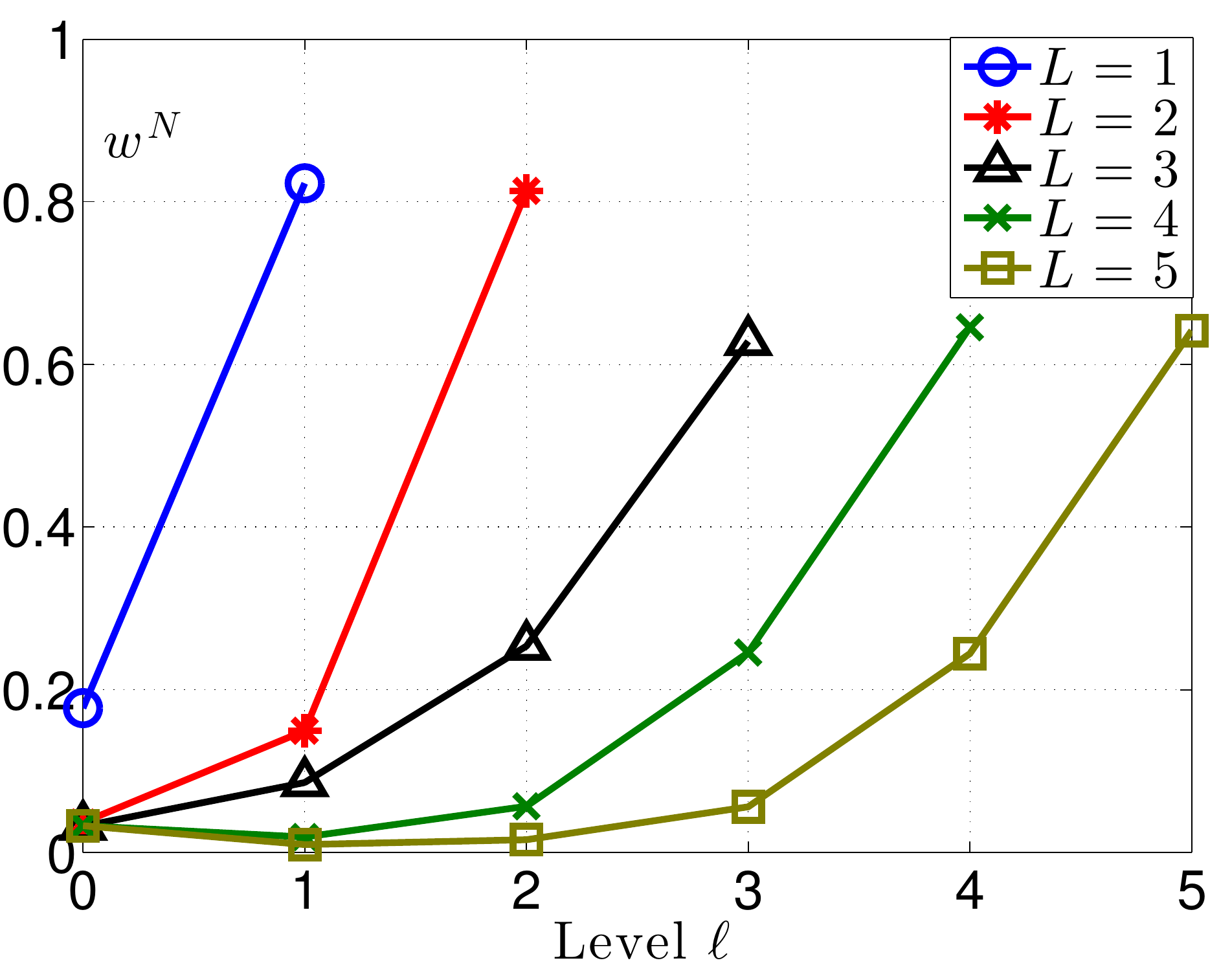}}
  \hfill \subcaptionbox{Expectation estimator with 95\% confidence interval.\label{fig:L-accuracyE}}
   [4.92cm]{\includegraphics[width = 4.85cm]{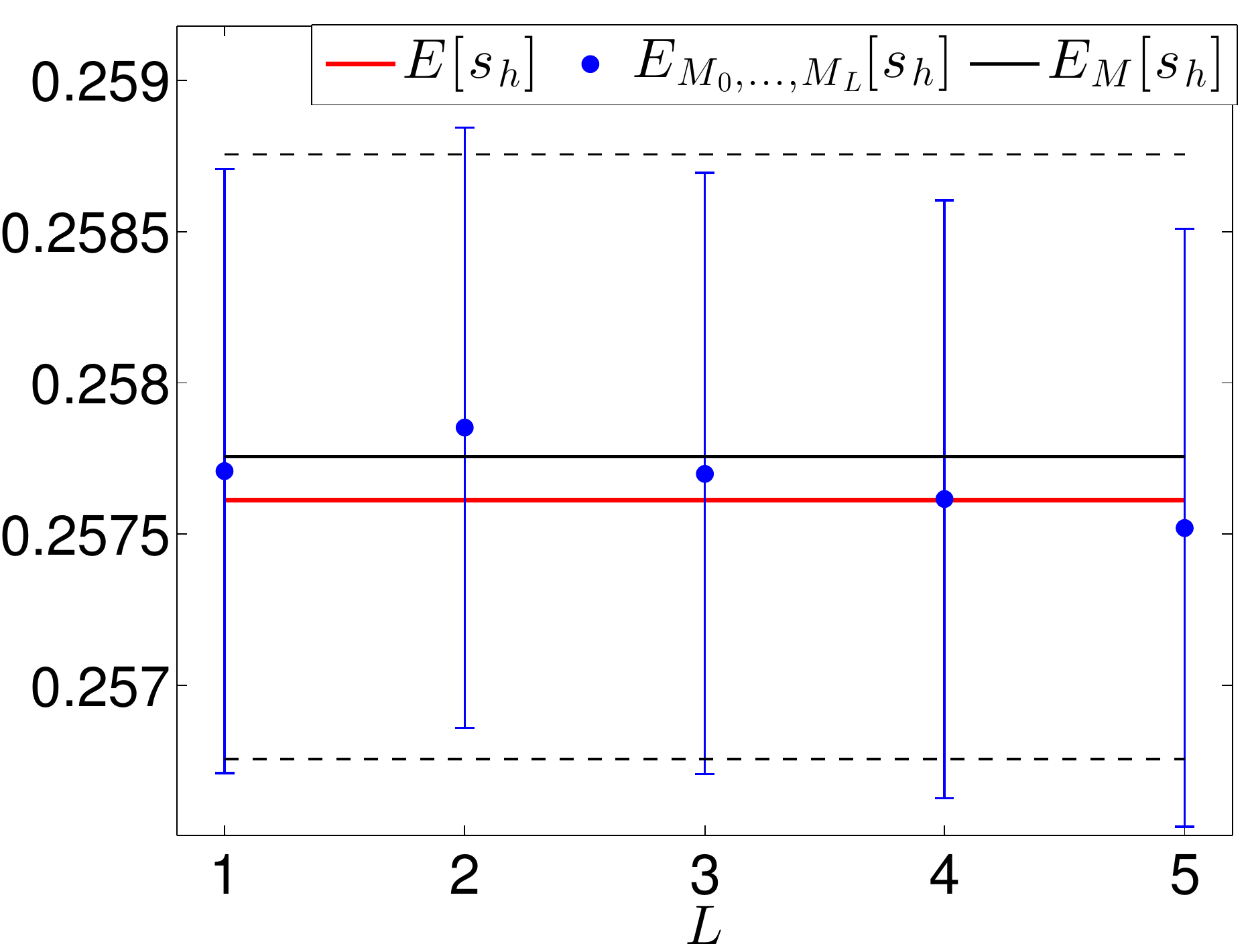}}
   \hfill \subcaptionbox{Variance estimator with 95\% confidence interval.\label{fig:L-accuracyV}}
    [4.92cm]{\includegraphics[width = 4.85cm]{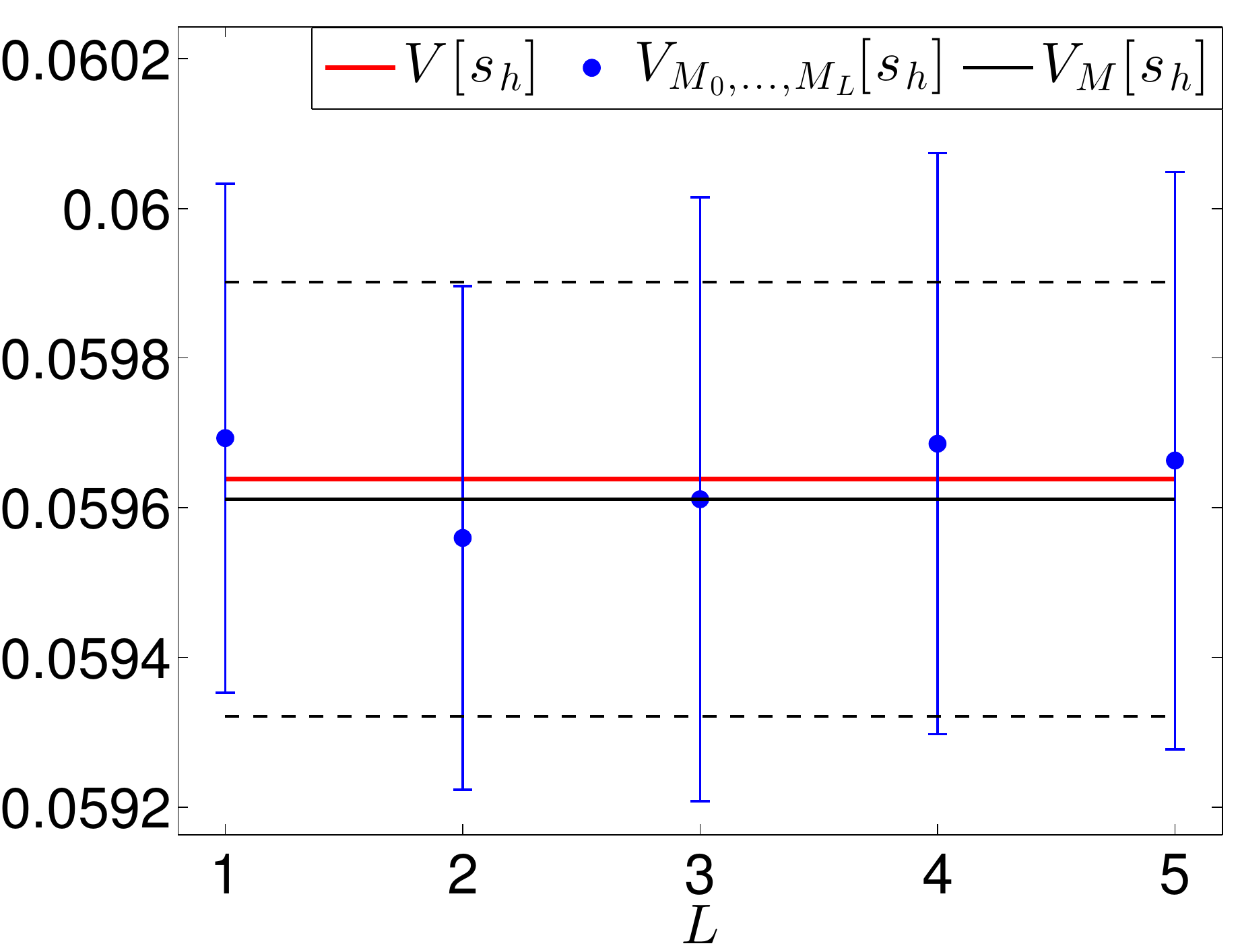}}
      \caption{Results for the $L$-MVR for arbitrary number of levels.}
    \end{figure}
We now analyze the performance of the $L$-MVR when an arbitrary number of levels are considered. The goal is to use the information on the test set to select not only the level sizes, but also the optimal multilevel model. For simplicity, we consider the same test set with $\widehat{M} = 100$ used for the previous cases, and obtain for each number of levels the optimal level sizes $\bm N = \lp N_1,\ldots,N_L \rp$, weights $\bm w^{\bm N}  = \lp w^{\bm N}_0,\ldots,w^{\bm N}_L \rp$ and the equivalent cost $\widehat{{C}}_L$ in \eqref{eq:cost3}. We then perform the $L$-MVR that provides the actual computational cost ${C}_L$ in \eqref{eq:cost2}. Results corresponding to $L = 1,2,3,4,5$ are presented in Table \ref{tab:multilevel}, for a confidence level of 0.95 and averaged 8 times. The costs are normalized with respect to the minimum costs, attained with $L=4$ for this problem. The actual speedup $\pi_L$ in \eqref{eq:cost} with respect to MC-HDG is also shown. 

The proposed selection method effectively predicts the \apr performance for each model, since $\widehat{C}_L/\widehat{C}_4$ replicates the behavior of ${C}_L/{C}_4$ quite well. The consideration of an arbitrary number of levels recovers even greater speedups, and the inexpensive \apr analysis enables the detection of the optimal model. Furthermore, the weights for each model are shown in Figure \ref{fig:L-weights}, exhibiting a nonlinear behavior that truly motivates its selection in an optimal automated manner. The distribution of the weights enforces a larger error on the coarser levels, and requires a smaller error on the finer levels relying on the reduction of variance. The estimators for each number of levels also satisfy the required accuracy, as seen in Figure \ref{fig:L-accuracyE}-\ref{fig:L-accuracyV}. 
   
\section{Conclusions}\label{sec:conclusions}

We have presented a model and variance reduction method for computing statistical outputs of stochastic elliptic PDEs. We first combined the reduced basis method with the hybridizable discontinuous Galerkin method by introducing a new HDG weak formulation that retains affine parametric dependence, hence providing rapid and accurate evaluation of the functional output of parametrized PDEs. We next incorporated them into the multilevel control variate framework to exploit the statistical correlation between the RB approximation and the high-fidelity HDG discretization to accelerate the convergence rate of the Monte Carlo simulations by several orders of magnitude. We then introduced {\em a posteriori} error bounds for the estimates of the statistical outputs. Finally, we devised an algorithm to select the RB dimensions and the number of levels $L$. We presented numerical results for both coercive and noncoercive elliptic problems. The results showed that the present method provides  a significant speedup compared  to both the 
MC-HDG method and the MC-RB method.

We conclude the paper by pointing out several possible extensions and directions for further research. Firstly, it would be interesting to address the computation of higher order moments, with the additional difficulties of determining the bias and the limiting distributions. Secondly, we would like to extend the proposed approach to nonlinear stochastic problems, which will broaden the application domain of our method. In this aspect, the main challenge remains the development of the RB method for the HDG discretization of nonlinear parametrized PDEs.  We would also like to tackle stochastic optimization problems with stochastic PDE constraints, for which the rapid and reliable evaluation of statistical outputs and their derivatives are crucial to finding an optimal solution of any stochastic optimization problem. We would like to develop new methods that allow us to compute not only the statistical outputs but also their derivatives with respect to the decision variables. 

\section*{Acknowledgements}
We would like to thank Dr. Xevi Roca and Professor Youssef Marzouk for countless fruitful conversations, suggestions and comments.

\appendix
\section{Bias of estimators}\label{sec:bias}
To simplify the notation, we will use the following auxiliary variables
\begin{alignat*}{2}
\zeta_h & := \lp s_h - E_{M_0,\ldots,M_L}[s_h] \rp^2,\qquad  &\zeta_{N_\ell} & := \lp s_{N_\ell} - E_{M_0,\ldots,M_L}[s_h] \rp^2,\; \ell = 1,\ldots,L\:,\\
\widehat{\zeta}_h & := \lp s_h - E[s_h] \rp^2,\qquad  &\widehat{\zeta}_{N_\ell} & := \lp s_{N_\ell} - E[s_h] \rp^2,\; \ell = 1,\ldots,L\:,\\
\overline{s}_h &= E_{M_0,\ldots,M_L}[s_h]\:,
\end{alignat*}
and the auxiliary (unbiased) variance
\beq
\widehat{V}_{M_0,\ldots,M_L}[s_h] = E_{M_0} [\widehat{\zeta}_h-\widehat{\zeta}_{N_1}] + \sum_{\ell = 1}^{L-1} E_{M_\ell} [\widehat{\zeta}_{N_\ell}-\widehat{\zeta}_{N_{\ell+1}}] + E_{M_L} [\widehat{\zeta}_{N_L}]
\eeq
which allows us to express the $L$-MVR variance estimate as
\beql\label{eq:apvar}
V_{M_0,\ldots,M_L}[s_h] = \widehat{V}_{M_0,\ldots,M_L}[s_h] -\lp E[s_h] - E_{M_0,\ldots,M_L}[s_h] \rp^2
\eeql
We first show the latter expression. We add and subtract $E[s_h]$ from every term within the expectations
\beq
\begin{split}
 V_{M_0,\ldots,M_L}[s_h] &= E_{M_0} \lb \blp s_h-E[s_h] - \overline{s}_h + E[s_h]\brp^2 -\blp s_{N_1}-E[s_h] - \overline{s}_h + E[s_h]\brp^2\rb\\
 &+ \sum_{\ell = 1}^{L-1} E_{M_\ell} \lb \blp s_{N_\ell}-E[s_h] - \overline{s}_h + E[s_h]\brp^2 -\blp s_{N_{\ell+1}}-E[s_h] - \overline{s}_h + E[s_h]\brp^2\rb\\
 &+ E_{M_L}\lb \blp s_{N_L}-E[s_h]- \overline{s}_h + E[s_h]\brp^2\rb \:,
\end{split}
\eeq
and now expanding the squares we arrive at
\beq
\begin{split}
 V_{M_0,\ldots,M_L}[s_h] &= E_{M_0} \lb \widehat{\zeta}_h - \widehat{\zeta}_{N_1} - 2\lp \overline{s}_h -E[s_h]\rp\lp s_h-s_{N_1}\rp\rb\\
 &+ \sum_{\ell = 1}^{L-1} E_{M_\ell} \lb \widehat{\zeta}_{N_\ell} - \widehat{\zeta}_{N_{\ell+1}} - 2\lp \overline{s}_h -E[s_h]\rp\lp s_{N_\ell}-s_{N_{\ell+1}}\rp\rb\\
 &+ E_{M_L}\lb \widehat{\zeta}_{N_L} - 2\lp \overline{s}_h -E[s_h]\rp\lp s_{N_L}-E[s_h]\rp + \lp\overline{s}_h -E[s_h] \rp^2\rb \:.
\end{split}
\eeq
Applying linearity of the MC expectation operator and grouping terms we arrive at \eqref{eq:apvar}. The bias of the $L$-MVR variance estimate is defined as
\beql\label{eq:ap1}
\begin{split}
E \lb V_{M_0,\ldots,M_L}[s_h] -V[s_h] \rb &= E \lb \widehat{V}_{M_0,\ldots,M_L}[s_h] -V[s_h] \rb - E\lb \lp E[s_h] - E_{M_0,\ldots,M_L}[s_h] \rp^2 \rb \\
&= E[s_h]^2 -E\lb E^2_{M_0,\ldots,M_L}[s_h] \rb. 
\end{split}
\eeql
since $E_{M_0,\ldots,M_L}[s_h]$ and $\widehat{V}_{M_0,\ldots,M_L}[s_h]$ are unbiased. If we rename the RB output differences as
\beq
z_0 := s_h - s_{N_1},\quad z_L:= s_{N_L},\quad z_\ell := s_{N_\ell} - s_{N_{\ell+1}},\; \ell = 1,\ldots,L-1
\eeq
the expression for expectation of the square of the $L$-MVR expectation estimate reads
\beq
  E\lb E^2_{M_0,\ldots,M_L}[s_h]\rb = \sum_{\ell = 0}^{L} \frac{1}{M^2_\ell}E\lb\lp \sum_{m=1}^{M_\ell} z_\ell(\sbf{y}_m)\rp^2\rb +2\sum_{\substack{\ell<\ell'\\ \ell =0}}^L E\lb\sum_{m=1}^{M_\ell}  z_\ell(\sbf{y}_m) \sum_{m'=1}^{M_{\ell'}} z_{\ell'}(\sbf{y}_{m'})\rb
\eeq
thanks to the linearity of the expectation operator. The latter expression can be further reduced with
\begin{align*}
 E\lb\lp \sum_{m=1}^{M_\ell} z_\ell(\bm y_m)\rp^2\rb &= M_\ell E[z_\ell^2] + (M_\ell^2 - M_\ell)E[z_\ell]^2,\\
 E\lb\sum_{m=1}^{M_\ell}  z_\ell(\sbf{y}_m) \sum_{m'=1}^{M_{\ell'}}z_{\ell'}(\sbf{y}_{m'})\rb  & = E[z_\ell]E[z_{\ell'}]\:,
\end{align*}
that hold because we consider independent samples within each level and independent samples among levels. We then have
\beq
  E\lb E^2_{M_0,\ldots,M_L}[s_h]\rb = \sum_{\ell = 0}^{L} \lp\frac{E[z_\ell^2] - E[z_\ell]^2}{M_\ell} + E[z_\ell]^2  \rp+2\sum_{\substack{\ell<\ell'\\ \ell =0}}^L E[z_\ell]E[z_{\ell'}]\:,
\eeq
and by induction on the number of levels, the cross-products can be reduced to
\beql\label{eq:ap2}
  E\lb E^2_{M_0,\ldots,M_L}[s_h]\rb = E[s_h]^2 + \sum_{\ell = 0}^{L} \frac{E[z_\ell^2] - E[z_\ell]^2}{M_\ell} = E[s_h]^2 + \sum_{\ell = 0}^{L} \frac{V[z_\ell]}{M_\ell}\:.
\eeql
Hence, if we combine equations \eqref{eq:ap1}--\eqref{eq:ap2} we obtain the bias of the $L$-MVR variance estimate
\beq
E \lb V_{M_0,\ldots,M_L}[s_h] -V[s_h] \rb =  -\sum_{\ell = 0}^{L} \frac{V[z_\ell]}{M_\ell}\:.
\eeq

\bibliography{mainbib}
\bibliographystyle{elsarticle-harv}

\end{document}